\documentclass[final]{siamart190516}
\usepackage[utf8]{inputenc}

\title{Kernel Analog Forecasting:\\ Multiscale Test Problems} 

\author{Dmitry Burov\thanks{Computing + Mathematical Sciences, California Institute of Technology, Pasadena, CA 91125 (\email{dburov@caltech.edu}, \email{kmanohar@caltech.edu}, \email{astuart@caltech.edu}).}
\and Dimitrios Giannakis\thanks{Department of Mathematics and Center for Atmosphere Ocean Science, Courant Institute of Mathematical Sciences, New York University, New York, NY 10012 (\email{dimitris@cims.nyu.edu}).}
\and Krithika Manohar\footnotemark[1]
\and Andrew Stuart\footnotemark[1]}

\usepackage{graphicx}
\usepackage{float}
\usepackage{overpic}
\usepackage{mathtools} 
\usepackage{amssymb}
\usepackage{tikz-cd}
\usepackage{xcolor}
\usepackage{bm}
\usepackage[shortlabels]{enumitem}
\usepackage{appendix}

\renewcommand\fbox{\fcolorbox{black}{white}}

\newcommand{\RMSE}{{\sf{{RMSE}}\,}}
\newcommand{\ac}{v} 
\newcommand{\hgp}{\ac_{\scriptscriptstyle{GP}}}
\newcommand{\cgp}{c_{\scriptscriptstyle{GP}}}
\newcommand{\ix}{{\color{black}\Pi}}
\newcommand{\est}[1]{{#1}_{\scriptscriptstyle{N}}}
\newcommand{\ltm}{L^2_{\mu}(\Omega;\R)}
\newcommand{\ltn}{L^2_{\mu_N}(\Omega;\R)}
\newcommand{\R}{\mathbb{R}}
\newcommand{\X}{\mathcal{X}}
\newcommand{\Y}{\mathcal{Y}}
\newcommand{\kerneleps}{\epsilon}
\newcommand{\lorenzeps}{\varepsilon}
\newcommand{\LNS}{L-96}

\DeclareMathOperator{\var}{var}
\DeclareMathOperator{\divr}{div}

\DeclareMathOperator{\spn}{span}

\DeclareMathOperator{\E}{\mathbb{E}}

\definecolor{lasunset4}{HTML}{6c5b7b}
\definecolor{lasunset5}{HTML}{355c7d}
\definecolor{egyptian}{HTML}{1443d4}
\definecolor{darkred}{rgb}{.7,0,0}

\definecolor{darkgreen}{rgb}{0,0.7,0}

\definecolor{darkblue}{rgb}{0,0,0.7}

\newcommand{\todo}[1]{{\color{lasunset4}{#1}}}
\newcommand{\revisionOne}[1]{{\color{black}{#1}}}

\theoremstyle{remark}

\newcommand{\Expect}{\operatorname{\mathbb{E}}}

\DeclareMathOperator*{\argmin}{arg\,min}

\DeclareMathOperator*{\proj}{proj}

\DeclareMathOperator*{\grad}{grad}

\newcommand{\bPhi}{{\Phi}}
\newcommand{\bPsi}{{\Psi}}
\newcommand{\bLambda}{{\Lambda}}
\newcommand{\ff}{{f}}
\newcommand{\bc}{{c}}

\newcommand{\xo}{\hat{x}}
\newcommand{\fo}{\hat{f}}
\newcommand{\bz}{{z}}
\newcommand{\bK}{{K}}
\newcommand{\ones}{\vec{1}}

\graphicspath{{./figures/}}

\begin{document}

\maketitle

\begin{abstract}
Data-driven prediction is becoming increasingly widespread as the
volume of data available grows and as algorithmic development matches
this growth. The nature of the predictions made, and the manner in which
they should be interpreted, depends crucially on the extent to which 
the variables chosen for prediction are Markovian, or approximately
Markovian. Multiscale systems provide a framework in which this issue 
can be analyzed. In this work kernel analog forecasting methods are studied 
from the perspective of data generated by multiscale dynamical systems.
The problems chosen exhibit a variety of different Markovian closures, 
using both averaging and homogenization; furthermore, settings where
scale-separation is not present and the predicted variables
are non-Markovian, are also considered. The studies provide guidance for the 
interpretation of data-driven prediction methods when used in practice.

\end{abstract}



\section{Introduction}

Data-driven prediction holds great promise in many areas of science
and engineering. Growth in the volume of data available in numerous
application areas has been matched by advances in computational methodologies
which are designed to utilize this data for prediction. However
fundamental questions arise in this field relating to the choice of
variables on which to base prediction, and whether or not the system
is Markovian in the chosen variables.  Whilst delay embedding can be
used to enhance the choice of variables in which Markovian structure
is present, prediction is often undertaken using variables in which
there is not a Markovian closure or in which this closure is only
approximate. The objective of the paper is to use multiscale systems 
to provide a framework in which the fundamental issue of the role
of Markovianity in data-driven prediction can be studied.  
We work within the setting of kernel analog
forecasting (KAF), a methodology that has seen success in
a number of application domains, and which is backed by a mature
theory. Subsection~\ref{ssec:LR} provides an overview of 
relevant literature in data-driven prediction for dynamical systems
and the multiscale setting in which we work. 
We outline our contributions to the understanding of
data-driven prediction within multiscale systems 
in Subsection~\ref{ssec:OC}.

\subsection{Background And Literature Review}
\label{ssec:LR}

In 1969, Lorenz originally introduced the idea of \emph{analog forecasting} for
prediction of dynamical systems using historical
data~\cite{lorenz1969atmospheric}.
Given initial data, the method locates its closest analog among the historical
points and reports the historical value of the corresponding observable, shifted
by the desired lead time.
By construction, analog forecasting avoids model error, but the resulting
forecast is not continuous with respect to initial data. It is therefore
non-physical and this fact, when combined with the paucity of data available 
at the time the method was proposed, limited the value of the methodology
in practice. An exponential growth in data volume has precipitated the
development of improved methodologies which build on Lorenz's original
idea, leading to algorithms which are backed by large data theories and
which depend smoothly on initial condition. This has been achieved through
the use of kernel based methods which result in data-driven prediction
based on weighting all the historical data according to its similarity 
to the initial data; this leads to algorithms which enforce continuity of the
forecast with respect to initial 
data~\cite{zhao2016analog} and, building on
the theory of reproducing kernel Hilbert spaces \cite{aronszajn1950theory}, 
to algorithms which can be theoretically justified in the large data 
limit \cite{Giannakis19,AlexanderGiannakis20}. 

KAF draws upon several fundamental ideas rooted in kernel methods for machine learning. 
First, the choice of kernel function is guided by the need for dimension reduction of big data and the specific learning task at hand. 
For clustering, similarity kernels~\cite{ScholkopfEtAl98} are used to construct graphs over the data, and clusters determined by its graph Laplacian eigenvectors~\cite{belkin2003laplacian}. 
This graph Laplacian construction is generalized in~\cite{coifman2006diffusion} to characterize diffusion operators on the manifold upon which the data lie via their eigenfunctions, known as diffusion maps, and further generalized in \cite{BerryHarlim16} to a class of variable-bandwidth kernels that control for variations in the density of the sampling distribution. Under different choices of kernel and normalization, the resulting eigenmaps can describe slow coordinates in dynamical systems~\cite{nadler2006diffusion_rc}, and can also be used as a basis to approximate evolution operators in stochastic differential
equations (SDEs)~\cite{berry2015nonparametric}. 
Algorithmic development has been aided by advances in the theory for 
pointwise~\cite{singer2006graph} and 
spectral~\cite{von2008consistency,TrillosEtAl19} 
convergence of these coordinates in the large data limit. 
Secondly, Markov operators constructed from kernels map into 
reproducing kernel Hilbert spaces (RKHSs)
in which kernel evaluation corresponds to function evaluation, and allow evaluating these eigencoordinates on out-of-sample data~\cite{coifman2006geometric} in a procedure known as Nystr\"om extension. This has found use in semi-supervised classification using support vector machines~\cite{scholkopf2002learning} to extend labels to new data, spline interpolation~\cite{wahba1990spline}, and forecasting~\cite{zhao2016analog}. 

In addition to exploiting ideas from kernel methods for machine learning,
 KAF may exploit additional structure from time-ordered data arising from
dynamical systems.  For example delay embedding is frequently used to 
identify Markovian structure. 
In diffusion forecasting~\cite{berry2015nonparametric}, diffusion maps are time-shifted to approximate the action of a shift operator on observables in 
SDEs. 
Alternatively, time-shifted diffusion maps can be used to approximate the reduction coordinates of this shift operator directly~\cite{Giannakis19}, or, temporal structure can be directly embedded into specialized cone kernels for analog forecasting~\cite{zhao2016analog}. 
The latter takes the Nystr\"om extension perspective of KAF, while in fact, KAF evaluates a conditional expectation of this shift operator, conditioned on the observations~\cite{AlexanderGiannakis20}.
The aforementioned shift operator known as the {\em Koopman} operator acts by composing observables with the dynamical flow map, and is a linear operator on these function spaces. 
Hence, the data-driven approximation of the Koopman operator is an exciting area of research.

Bernhard Koopman introduced the linear operator that carries his 
name in the 1930s as part of his study of ergodic and 
Hamiltonian dynamics~\cite{koopman1931hamiltonian}.  
In data-driven identification of coherent structures 
spectral decompositions of the Koopman operator \cite{Mezic05} and 
the related transfer operator \cite{DellnitzJunge99} play a central
role, driven by the fact that in such a basis, forecasting of nonlinear dynamics amounts to scalar multiplication by eigenvalues. 
Algorithms used in practice compute finite-dimensional regression onto pre-computed libraries consisting of functions of the time-lagged snapshots~\cite{schmid2010dynamic,rowley2009spectral,klus2018data}. 
However, convergence guarantees are limited, requiring
stringent assumptions on the libraries and spectrum~\cite{arbabi2017ergodic}. 
In particular mixed spectra resulting from chaotic/mixing 
systems pose a challenge for numerical methods. 
Recent data-driven methods which leverage infinite-dimensional 
feature spaces provided by kernels \cite{giannakis2018reproducing}, 
as well as kernel constructions in spectral space \cite{KordaEtAl18}, 
are able to tackle the continuous 
part of the spectrum of the Koopman operator. 
For forecasting purposes, pointwise evaluation of the Koopman operator 
acting on observables is the natural setting, rather than spectral
approximation, and is the perspective we take.

%
%
%

Multiscale analysis provides a setting in which to understand the role
of rapidly varying (in space or time)
system components on the slowly varying
variables used for predictive models
\cite{weinan2011principles}.
In this paper we will work in the framework of averaging and homogenization for
partial differential equations (PDEs) and SDEs, as developed
in~\cite{bensoussan2011asymptotic}.
\revisionOne{Chapters 9, 10 and 11 of the 
book \cite{pavliotis2008multiscale} contain a
pedagogical exposition of the subject that is adapted to the chaotic 
deterministic ordinary differential equations 
setting that is the focus of this paper.
However the rigorous extension of the theory of averaging and 
homogenization to ODEs, rather than SDEs, is non-trivial and 
less well-developed. Early work in this direction was 
contained in~\cite{papanicolaou1974asymptotic}. However it was not 
until the fundamental work of Melbourne and
co-workers that a theoretical approach with verifiable conditions 
was developed 
\cite{melbourne2005almost,melbourne2011note,kelly2016smooth,kelly2017deterministic}.

We will use the example developed in \cite{melbourne2011note}, 
which exploits the (proven) chaotic properties of the Lorenz 63 
model \cite{lorenz1963deterministic,tucker1999lorenz}, 
to provide an example of a chaotic ordinary differential equation (ODE) 
which homogenizes to give an SDE.  
And we will use the multiscale Lorenz 96 model \cite{lorenz1996predictability}
to provide an example of an ODE to which the averaging principle
may be applied to effect dimension reduction, as pioneered and exploited
in \cite{fatkullin2004computational}; we note, however, that the mixing
properties required to \emph{prove} the averaging principle for
the Lorenz 96 model have not been established, even though numerical
evidence strongly suggests that it applies in certain parameter regimes.} 
The work of Jiang and Harlim \cite{jiang2019modeling} studies data-informed
model-driven prediction in partially observed ODEs, using
ideas from kernel based approximation; and in the paper \cite{harlim2019machine}
the idea is generalized to discrete time dynamical systems, and
neural networks and LSTM modeling is used in place of kernel methods. 
In both the papers \cite{jiang2019modeling,harlim2019machine} 
multiscale systems are used to test their methods in certain regimes.

Data-driven analog forecasting, kernel methods, and Koopman methodologies have each individually found widespread use in real-world forecasting and coherent pattern extraction applications. 
Analog forecasting, albeit without kernels, has been used to predict weather patterns~\cite{chattopadhyay2020analog,delle2013probabilistic,van1989new}, yet is known to have limitations predicting chaotic behavior. 
Khodkar et al. recently developed a Koopman-based framework 
using delay embedded observables to predict chaotic dynamics~\cite{khodkar2019koopman}. 
Nonlinear Laplacian spectral analysis, which applied kernel and delay embeddings akin to Koopman observables, successfully recovers \revisionOne{coherent oscillatory phenomena such as the seasonal/diurnal cycles and the El Ni\~no Southern oscillation~\cite{SlawinskaGiannakis17} through kernel eigenfunctions \cite{DasGiannakis19}. Implemented with such kernels, KAF would naturally capture the fundamental oscillatory components of the predictand variable, and thus interpolate between an initial-value (``weather'') forecast at short times and a climatology forecast at asymptotic times when the mixing component of the predictand has decayed; e.g.~\cite{WangEtAl20}.} 
Koopman operator approximation has also been widely adopted to study high-dimensional complex, even turbulent fluid flows \cite{GiannakisEtAl18}, see~\cite{mezic2013analysis} for a study of these applications.

In this paper we use KAF as developed in the papers 
\cite{AlexanderGiannakis20,Giannakis19,zhao2016analog}.
Technical details on the implementation, including pseudocode
and further relevant references, are collected 
in Appendix~\ref{appComputation}.

\subsection{Our Contribution}
\label{ssec:OC} 

We use multiscale methodology to introduce four classes of ODE
test problems which exhibit Markovian dynamics after
elimination of fast variables; stochastic, chaotic, quasiperiodic and
periodic behavior may be obtained in the slow variable, depending on the
setting considered. Using these test problems, our four
main contributions to the understanding of KAF are as follows:

\begin{enumerate}

\item We apply KAF techniques to data generated by each of these four classes
of multiscale test problems, and use the behaviour of the averaged or
homogenized slow system to interpret the resulting predictions. 
In particular, KAF methods converge, in the large data limit,
to a conditional expectation defined via the Koopman operator of
the multiscale systems; we use this as the basis for our interpretation. 
Moreover, we demonstrate the use of KAF to estimate the variance 
of the predictor itself. In each of the four cases the 2$\sigma$-interval 
captures the real trajectory,
even when the KAF predictor does not track the trajectory itself.
This can be viewed as a separate application of the KAF methodology 
which will be useful in cases when forecasting of high probability
bounding sets suffices even when the trajectory itself is hard to predict. 
In all cases we also study problems in which the scale-separation
is not present, but KAF prediction of mean and variance is attempted 
on the basis of data from only a subset of the variables.

\item The KAF method is based on data-driven approximation of the 
eigenvalue problem arising from a kernel integral operator. In the 
setting in which the multiscale ODE homogenizes to produce an SDE 
corresponding to a bistable gradient system with additive noise,
a limiting analytical expression is available for the eigenfunctions;
we demonstrate that these limiting eigenfunctions are well-approximated
by the data-driven method. This comparison gives insight into the 
empirical methods used to tune free parameters within KAF. 

\item  In the setting in which the multiscale ODE averages to produce 
an ODE of lower dimension we use alternative data-driven ODE closures as
a benchmark against which to compare the purely data-driven KAF methods. 
This gives insight into the relative merits of purely data-driven
prediction, and prediction which combines model-based knowledge with
data.

\item We use the insights from these carefully constructed numerical
experiments to make recommendations about deployment, and parameter-tuning,
of KAF methods to real data.

\end{enumerate}

The paper is organized as follows. In Section~\ref{sec:M}, we outline 
the data-driven construction of the prediction function using KAF methodology.
We explain the sense in which the construction converges to a conditional 
expectation defined via the Koopman operator associated to a measure-preserving
dynamical system assumed to underlie the data. We also describe two kinds of 
canonical multiscale systems which give rise to homogenization and averaging
effects, and which we use to provide interpretation of this conditional 
expectation.  Section~\ref{sec:63} introduces a test problem in the form
of a double-well gradient flow driven by chaotic Lorenz 63 dynamics which 
homogenizes to give an SDE in the scale-separated regime; numerical results 
applied to prediction of the slow variable exhibit contributions 1 and 2. 
In section~\ref{sec:96}, we introduce the  multiscale Lorenz 96 system 
which averages to give an ODE in the slow variables;  three different 
parametric regimes give rise to periodic, quasiperiodic, and chaotic responses
in the slow variable.  The behaviour of KAF-based prediction in these three
regimes is studied, to illustrate contribution 1; and a slow-variable
closure model, built using Gaussian process regression, is compared with
the KAF to illustrate contribution 3. In section \ref{sec:C} we overview
the insights obtained by studying KAF methods through the lens of
multiscale systems; and we then make concrete recommendations about 
interpreting the output of KAF techniques when applied to naturally
occurring data, contribution 4.

\section{Methodology} 
\label{sec:M}

In subsection \ref{ssec:OOM} we overview the two key ideas which
interact to underpin the studies in this paper: KAF and multiscale methods, 
tailoring the exposition to the use of the latter as a tool to
understand the former. 
We then give more details on KAF.
The two primary components of the KAF methodology are: 
(i) viewing forecasting as evaluation of a conditional expectation 
of the Koopman operator applied to the desired observable; 
(ii) approximation of this
conditional expectation in a data-driven fashion. Subsections \ref{ssec:KF}
and \ref{ssec:DDA} describe (i) and (ii) respectively, whilst
subsection \ref{as:kernel} is devoted to a key practical component
of the data-driven approximation, namely construction of the kernel,
and subsection \ref{as:L} to the data-driven choice of integer $\ell$,
the number of (approximate) eigenfunctions used in the data-driven forecast.

\subsection{Overview Of Methodology}
\label{ssec:OOM}

The problem setting for prediction is as follows. 
We assume that we are given $N$ time-ordered data samples  
\begin{displaymath}
\{x_n\}_{n=0}^{N-1} \subset \mathcal{X},
\end{displaymath}
where $x: \R \to \mathcal{X}$ 
is a continuous time process,  $x_n=x(n\Delta t)$ and $\Delta t$ 
is the sampling rate. We assume that the continuous time 
process $x$ in $\mathcal{X}$ is derived from Markovian dynamics for
a coupled pair $(x,y)$ evolving in  
the larger state space $\mathcal{X}\times\mathcal{Y}$. 
Assume that the desired prediction {\em lead time} $\tau$ is an integer multiple of the sampling interval, that is, $\tau = q\Delta t$. Included with the data are values of the associated 
{\em prediction observable} advanced by $\tau$ time units 
\begin{displaymath}
\{f_{n+q}\}_{n=0}^{N-1} \subset \R,
\end{displaymath}
defined by the Markovian dynamics via an 
unknown map $F\colon \X\times\Y \to\R$; thus
$f_n = F(x_n,y_n).$
The goal of KAF is to predict $F(x(\tau),y(\tau))$ given only partial information, $x(0) = x$, and  
the $N$ data samples $x_n$. We view the data-driven predictor as
a map $Z_{\tau}\colon \X \to \R$ which takes initial condition
$x$ as input.

Given initial data $x$ and lead time $\tau$, the KAF predictor averages over the $\tau$-shifted observable values provided in the training data and weighted by a kernel $p:\mathcal{X}\times\mathcal{X}\rightarrow \R$ constructed from
the data; the resulting algorithm has the following form:
\begin{equation}
\label{eq:Z_EXPAND}
\begin{aligned}
Z_{\tau}(x) &= \frac{1}{N}\sum_{n=0}^{N-1} p(x,x_n) f_{n+q},\\
p(x,x_n)&=\sum_{j=0}^{\ell(\tau)-1} \frac{ \psi_j(x)\phi_j(x_n)}{\lambda_j^{1/2}}.
\end{aligned}
\end{equation}
The weighting kernel $p(x,x_n)$ determines how much weight to
attach to a time-series initialized at point $x_n$, according to
its proximity to $x$, the desired initial point.
The features $\phi_j$ are computed from an eigenvalue problem 
associated with a data-driven
approximation of a kernel integral operator, constructed from $x_n$;
in the large data limit this provides an orthonormal basis for the
entire space. The function $\psi_j$ is an out-of-sample Nystr\"om  
extension of $\phi_j$, orthonormalized with respect 
to an underlying RKHS structure. 
The method may be seen as a smoothed version of Lorenz's original proposal 
for data-driven prediction -- analog forecasting \cite{lorenz1969atmospheric}. \revisionOne{Analog forecasting, by contrast, predicts the trajectory in the training data obtained by finding the training data point nearest to the given initial condition in some metric $d:\X\times\X\rightarrow\R$
\begin{equation}
\label{eq:LAF}
Z_{\tau}(x) = f_{n^{\star} + q}, \quad \quad n^\star = \argmin_{n=0,\dots,N-1} d(x,x_n).
\end{equation}
It can be seen that Lorenz' method will result in predictions discontinuous with respect to initial data, especially for systems that exhibit sensitive dependence on initial conditions. }
In particular, KAF addresses the issue of 
continuity of the prediction with respect to the initial condition, and 
it does so in a framework which is provably statistically consistent in 
the large data limit \cite{AlexanderGiannakis20,Giannakis19,zhao2016analog}. 
Further details of the methodology are given in the next 
two subsections, and the attendant information in Appendix \ref{a:a}.

An important challenge addressed by this methodology is that, since 
the $y$ component of the system is not observed, the sequences $\{x_n\}$ 
and $\{f_{n+q}\}$ are non-Markovian. As a consequence the 
standard idea of constructing a Markov chain from the data is not natural. 
The kernel analog forecasting method evaluates a conditional expectation 
of the forecast conditioned, using the observed data $\{x_n\}$, explicitly 
incorporating information loss resulting from unobserved $y$;
it is hence a natural approach to the problem at hand. 
Multiscale systems provide a natural setting for the study of KAF
methods, and in particular the issue of prediction of non-Markovian
or approximately Markovian systems.
In this paper we will consider the variable $x$ as the slow component
of a Markovian system for pair $(x,y)$ in which $y$ evolves
as a fast variable.  We consider averaging and homogenization
settings in which the dynamics 
for $x$ is approximately Markovian, and the conditional expectation 
arising in the KAF method may be understood explicitly. 
This will enable us to obtain a deeper understanding of how KAF
works, and help users of the methodology interpret it.
We now outline the averaging and homogenization settings that we will use. 
Details of the theory underlying them
may be found in~\cite{pavliotis2008multiscale}.

We will study multiscale systems which exhibit averaging, in the form
\begin{equation}
\label{eq:avg}
\begin{cases}
& \dot{x} = v_0(x) + By, \\
& \dot{y} = \frac{1}{\varepsilon}g(x,y),
\end{cases}.\tag{A}
\end{equation}
where \revisionOne{$B\colon \Y \to \X$} is linear. 
The average of $By$ under the invariant measure
of the $y$ dynamics, with $x$ frozen, provides a closed
approximate ODE dynamics for $x$, when $\varepsilon$ is small.
If we denote the $x-$parameterized invariant measure for the
$y$ dynamics with $x$ frozen by
$\nu^x(dy)$ then for $\lorenzeps \ll 1$ we obtain
$x \approx X$ where
\begin{equation}
  \label{eq:averaged_ode}
  \begin{aligned}
    \dot{X} &= v_0(X) + \ac(X), \\
    \ac(\zeta)    &= \int_\Y By \, \nu^\zeta(dy).
  \end{aligned}
\tag{A0}
\end{equation}
\revisionOne{To guarantee uniqueness of solutions in~\eqref{eq:averaged_ode}, it
suffices that the conditional measure has enough continuity, as a function of
$x$, so that $\ac(\zeta)$ is Lipschitz.}

When the variable $By$ averages to zero a different scaling is required
to elicit the effect of the fast variable on the slow one. To this end
we also consider multiscale systems which exhibit
homogenization, in the form 
\begin{equation}
\label{eq:homog}
\begin{cases}
& \dot{x} = v_0(x) + \frac{1}{\varepsilon}By \\
& \dot{y} = \frac{1}{\varepsilon^2}g(y)
\end{cases}\tag{H}
\end{equation}
Here we assume that 
$$\int_{\mathcal Y} By \, \nu(dy)=0$$ where $\nu$ is the 
invariant measure of the $y$ dynamics.
The approximate dynamics for $x$, when $\varepsilon$ is small,
is then governed by an SDE in this setting;
\revisionOne{the work of Melbourne and co-workers provides the sharpest results
in this context
\cite{melbourne2005almost,melbourne2011note,kelly2016smooth,kelly2017deterministic}}.
If this is the case then, invoking the homogenization 
principle, $x \approx X$ where $X$ is governed by an SDE of the form 
\begin{equation}
  \dot{X} = v_0(X) + \sqrt{2\sigma} \dot{W},
\tag{H0}
\end{equation}
where $W$ denotes the Wiener process, and $\sigma$ is a uniquely
determined positive constant that can be computed numerically from
the mixing properties of the $y$ process.

\subsection{Koopman Formulation Of Prediction}
\label{ssec:KF}

\begin{figure}[!ht]
  \label{fig:comm_diag}
  \centering
  \begin{minipage}[b]{0.4\textwidth}
    \centering
    \begin{tikzcd}[row sep=scriptsize,column sep=scriptsize]
      \X \times \Y \arrow[rr, "\Phi^\tau"] \arrow[ddrr, "U^\tau F"] & & \X \times \Y \arrow[dd, "F"] \\
      \\
      & & \R \\
    \end{tikzcd}

    Koopman
  \end{minipage}
  \begin{minipage}[b]{0.4\textwidth}
    \centering
    \begin{tikzcd}[row sep=scriptsize,column sep=scriptsize]
      \X \times \Y \arrow[rr, "H"] & & \X \arrow[ddll, bend left, "Z_\tau"] \\
      \\
      \R \\
    \end{tikzcd}

    Data-driven
  \end{minipage}

  \caption{
    Comparison of exact Koopman picture and data-driven approximation.
    The approximate mapping $Z_\tau(\cdot)$ 
is constructed from the data streams
$\{H\Phi^t(x_0,y_0)\}_{0 \leqslant t \leqslant T}$ and
$\{F(\Phi^t(x_0,y_0))\}_{0 \leqslant t \leqslant T}.$
  }
\end{figure}
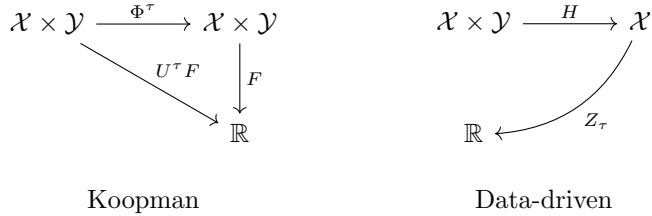

We let $\Omega= \X \times \Y$ and assume that $\Phi^t\colon \Omega \to \Omega$ 
is an ergodic dynamical system with invariant probability measure $\mu$; we 
assume $t \in \R^+$ but the extension to discrete time is straightforward. 
Define the continuous observation map $\ix \colon \Omega \to \X$ and the prediction
observable $F\colon \Omega \to \R$; we assume that $F$ is square-integrable with
respect to the invariant measure:
\[
  F\in L^2_{\mu}(\Omega;\R) \coloneqq \left\{ F\colon \Omega \to \R \; \middle|
  \; \int_{\Omega} |F(\omega)|^2 \mu(d\omega)<\infty \right\}.
\]
We define the Koopman operator $U^t \colon \ltm \to \ltm$ by $U^t g = g \circ
\Phi^t$.
We seek, in a sense to be made precise, the function $Z_\tau \colon \X \to \R$
such that $Z_\tau \circ \ix$ 
is the best approximation to a perfect prediction of
$F \circ \Phi^\tau$.
We formalize this by introducing the Hilbert subspace $V \subseteq \ltm$ given
by 
\[
  V \coloneqq \left\{ g \in \ltm \; \middle| \; g = g' \circ \ix, \; g'\colon \X \to \R \right\}.
\]
This Hilbert space captures the notion of making predictions based only
on information in $\X$. 
Note that the perfect forecast would satisfy $Z_\tau \circ \ix= U^\tau F$, but
that such a forecast will not be possible in general because $\X$ is a proper
subset of $\Omega$, and information needed for perfect prediction of $F$ will be
missing.
Among all elements of $V$, the minimal prediction error in $\ltm$ is attained
by the conditional expectation
\begin{gather}
  \label{eq:g}
  \E (U^\tau F \mid \ix) = \argmin_{g'\in V} \lVert g' - U^\tau F \rVert_{L^2_\mu} = \proj_{V} U^\tau F.
\end{gather} 
This formulation of prediction encapsulates the inherent loss of information
incurred through observing only a set of functionals of an ergodic
dynamical system, and the effect of this loss of information on prediction.  
In subsequent sections of this paper we will assume that 
$F = F' \circ \ix $ for some $F': \X \to \mathbb R$ because it is often
natural to try to predict only functionals of the slow variables. Note,
however, that the methodology is not restricted to such $F$ and in this 
subsection, the next subsection and in Appendix \ref{a:a}
we describe the more general setting for completeness.

\subsection{Data-Driven Approximation}
\label{ssec:DDA}

The formulation in the preceding subsection
encapsulates the inherent loss of predictive power incurred through 
observing only a set of functionals of an ergodic dynamical system.  
This is formalized by seeking the best approximation of the
Koopman evolution from within a Hilbert subspace capturing the
notion of depending only on specified functionals on the state
space of the dynamical system. We now demonstrate how data may be
used to further approximate this best approximation, and to do so
in a manner which is refineable as more data is acquired.
The approach is summarized in Figure \ref{fig:comm_diag}.

For observation time $\Delta t>0$ we define
\begin{align*}
\omega_n &= \Phi^{t_n}( \omega_0 ), \quad t_n = n\Delta t,\\
x_n &= \ix( \omega_n ), \quad f_n = F( \omega_n ).
\end{align*}
We assume that we are given time-ordered pairs
\begin{equation}
\label{eq:nd}
\{ ( x_0, f_q ), ( x_1, f_{1+q} ), \ldots, ( x_{N-1},  f_{N-1+q} ) \},
\end{equation}
and the objective is to construct, from this data, a function
$Z_\tau \colon \X \to \R$ which predicts $F$ at lead time $\tau$
so that $Z_\tau \circ \ix \approx g^\star$ where $g^\star$ solves the
minimization problem in \eqref{eq:g}.
Furthermore we wish to carry this out in a manner which ensures that, in an
appropriate topology, $Z_\tau \circ \ix \to g^\star$ as $N \to \infty$.

To this end we introduce a hypothesis space $\mathcal{H}_{\ell,N}$,
of dimension $\ell$ and depending on the $N-$dependent data set \eqref{eq:nd},
and seek to solve the minimization problem 
\begin{equation}
  \label{eq:LN}
  Z_\tau = \argmin_{g \in \mathcal{H}_{\ell,N}} \lVert g \circ \ix - U^\tau F \rVert_{L^2_\mu}. 
\end{equation}
The choice of the hypothesis space is constrained by the need to be
able to solve the minimization problem \eqref{eq:LN} explicitly, 
using only the data \eqref{eq:nd}, and by the requirement that
$Z_\tau \circ \ix$ recovers $g^\star$ in the large data limit $N \to \infty$. 
Moreover, in order to be practically useful forecast functions, elements of $\mathcal H_{\ell, N}$ should allow pointwise evaluation at any $ x \in \mathcal X$, which is not defined in arbitrary subspaces of $L^2_\mu(\Omega;\mathbb R)$.

With these considerations in mind, we introduce a kernel function $k\colon \X \times \X \to \R$
and RKHS $\mathcal{K}$ with the properties
\[
  f( x ) = \langle k_x, f \rangle_\mathcal{K}, \quad
  k_x = k( x, \cdot ), \quad
  \langle k_x, k_{x'} \rangle_{\mathcal{K}} = k(x,x').
\]
We then define $\mathcal H_{\ell,N}$ as an $\ell$-dimensional subspace of $\mathcal K$, to be described below. We also note that the kernel $k$
is constructed from a data-stream of length $N$, but we suppress the
explicit dependence of $k$ on $N$ in the notation. In Appendix \ref{a:a} 
we discuss our data-driven construction of $k$, and choice of $\ell$. 
For now we proceed on the assumption that we have a kernel, and hence 
a RKHS, as well as a method for choosing $\ell.$ 

Let $\est{\mu} = \frac{1}{N} \sum_{n=0}^{N-1} \delta_{\omega_n}$ be the sampling measure 
underlying the training data~\eqref{eq:nd} and define
\[
  \ltn \coloneqq \left\{F\colon \Omega \to \R \; \middle| \; \int_{\Omega} |F(\omega)|^2 \est{\mu}(d\omega)
  =\frac{1}{N}\sum |F(\omega_n)|^2<\infty \right\}.
\]
Associated with $\est{\mu}$ is an integral operator $G\colon \ltn \to \ltn$
which we identify with  a symmetric, positive-semidefinite, $N \times N$ kernel
matrix $G \in \R^{N \times N}$ with entries
\begin{equation}
\label{eq:G}
  G_{mn} =  k( x_m, x_n ), \quad x_n = \ix(\omega_n), 
\quad 0 \leq m,n \leq N -1.
\end{equation}
The eigenvectors of this matrix lead to an orthonormal basis 
$\{ \phi_j \}_{j=0}^{N-1} $ of $\R^N$ such that
\[
    G \phi_j = \lambda_j \phi_j, \quad \lambda_0 \geq \lambda_1 \geq \cdots \geq \lambda_{N-1}, \quad \lVert \phi_j \rVert_2 = \sqrt{N}.
\]
We may also identify each element $\phi_j \in \R^N$ with element
$\phi_j \circ \ix\in \ltn$ via the definition $\phi_j(x_n)$ as the $n^{th}$ entry
of the vector $\phi_j \in \R^N$. Using the
same symbols for elements of $\R^N$ and $\ltn$, as well as for
linear transformations on those spaces, is a useful economy of notation.
Then the following functions $\psi_j \colon \X \to \R $ 
form an orthonormal set in $\mathcal{K}$:
\begin{equation}
\label{eq:psiN}
  \psi_j = \frac{1}{N \lambda_j^{1/2}} \sum_{n=0}^{N-1} k( \cdot, x_n )\phi_j(x_n), \quad \lambda_j > 0.
\end{equation}
This is a form of Nystr\"om extension \cite{coifman2006geometric}.

As hypothesis space we take 
\begin{equation}
\label{eq:HS}
  \mathcal{H}_{\ell,N}=\spn\{\psi_0, \cdots, \psi_{\ell-1}\} \subseteq \mathcal{K}
\end{equation}
noting that the basis functions themselves depend on the data set, and hence on
$N$.
We may now solve the optimization problem \eqref{eq:LN} and an explicit
computation yields, for $\tau=qT$,
\begin{equation}
\label{eq:Z_PRED}
                Z_\tau(x) = \sum_{j=0}^{\ell-1} \frac{c_j(\tau)}{\lambda_j^{1/2}} \psi_j(x), \quad c_j(\tau) = \langle \phi_j \circ \ix, U^\tau F \rangle_{L^2(\est{\mu})} = \frac{1}{N} \sum_{n=0}^{N-1} \phi_j(x_n) f_{n+q}.
\end{equation}
Note that this construction of the predictor $Z_\tau$ is entirely data-driven:
the basis functions $\psi_j$ and the eigenvalues $\lambda_j$ are found
from the eigenvalues and eigenvectors
of the data-defined kernel matrix; and the coefficients $c_j$ are
computed as sums over the data set. Furthermore, Theorem 14 in 
\cite{AlexanderGiannakis20} proves that $Z_\tau \circ X$ converges to 
$g^\star$, the solution of the minimization problem \eqref{eq:g}, 
as $N \to \infty$, followed by $\ell \to \infty$, in an $L^2$ sense 
with respect to the invariant measure $\mu$ on $\Omega.$

More generally, any function of the observable can be predicted in this 
data-driven manner, which provides a convenient framework for uncertainty 
quantification. The conditional variance 
between forecast and ground truth can also be computed in the hypothesis space as in~\eqref{eq:Z_PRED} using the coefficients
\begin{equation}
\label{eq:V_PRED}
             \hat{c}_j(\tau) = \langle \phi_j \circ \ix, (U^\tau F-Z_\tau)^2 \rangle_{L^2(\est{\mu})} = \frac{1}{N} \sum_{n=0}^{N-1} \phi_j(x_n)(f_{n+q}-Z_\tau(x_n))^2.
\end{equation}
For detail on the data-driven kernel construction, the data-driven choice
of $\ell$ and the conditional variance estimator see subsections \ref{as:kernel},
\ref{as:L} and \ref{as:V} respectively.

%

\section{Homogenization: Lorenz 63 Driven System}
\label{sec:63}

This section is devoted to the setting in which a chaotic ODE
of form (H) is approximated by an SDE of form (H0).
The goal is to make predictions of the $x$ variable, using 
data concerning only the $x$ variable from (H); 
the role of (H0) is simply to help us interpret those
predictions. This setting presents unique challenges for forecasting as 
one cannot expect the outcome of any method to predict a sample path of a
stochastic process without knowledge of the driving noise. This fact
has direct bearing on prediction in (H) using $x$-data alone, since
(H0) demonstrates that the time series of the $x-$data is approximately 
that of an SDE; without knowledge of the noise, which is governed by
the unobserved $y$ variable, prediction of the trajectory of $x$ is not
possible. 
In subsection~\ref{ssec:CE63} we examine instead the long-term statistics
predicted by KAF from data generated by (H) --- the conditional expectation 
and variance of the stochastic process --- and compare them with estimates 
computed from (H0) using Monte-Carlo simulation of the SDE. This illustrates
our main contribution 1 from the list in subsection \ref{ssec:OC}. 
Then, in subsection~\ref{ssec:E}, exploiting the fact that the limiting 
process is one-dimensional, we find explicit expressions for the 
kernel eigenfunctions in the limit problem (H0) and compare these with the
eigenfunctions obtained from data-driven techniques applied to (H), 
our main contribution 2 from subsection~\ref{ssec:OC}.
Subsection~\ref{ssec:NM1} is concerned with non-Markovian prediction,
in which there is no scale-separation between observed and unobserved
variables. 
We start, however, in subsection~\ref{ssec:L63}, introducing the
concrete model around which our experiments are organized.

\subsection{The Model}
\label{ssec:L63}

The first test problem arises from driving a double-well gradient
flow with a chaotic signal obtained from the \revisionOne{Lorenz 63 model
\cite{givon2004extracting}}:
\begin{equation}
\label{eq:l63}
\begin{aligned}
    \dot{x} &= x-x^3 + \frac{4}{90\lorenzeps}y_2,  \\
    \dot{y}_1 &= \frac{10}{\lorenzeps^2}(y_2-y_1), \\
    \dot{y}_2 &= \frac{1}{\lorenzeps^2}(28y_1-y_2-y_1y_3), \\
    \dot{y}_3 &= \frac{1}{\lorenzeps^2}\left(y_1y_2- \frac{8}{3}y_3\right).
\end{aligned}
\end{equation}
This is of form (H).
In \cite{melbourne2011note} it is proved that as $\lorenzeps \to 0$ this
system converges weakly in $C([0,T];\mathbb{R})$, when projected onto the $x$ variable alone,
to solution of the SDE 
\begin{equation}
\label{eq:SDE}
\begin{aligned}
\dot{X} &= -\Xi'(X) + \sqrt{2\sigma}\dot{W},\\
    \Xi(X) &= \frac{1}{4}(X-X^2)^2.
\end{aligned}
\end{equation}
Thus, this white-noise driven gradient system is of form (H0).
The value of the constant $\sigma$ is identified in \cite{givon2004extracting}.
For the current work
the key point to appreciate is that for small $\lorenzeps$ the variable $x$
in \eqref{eq:l63} exhibits (approximately) Markovian behaviour, 
but this behaviour is stochastic. The SDE is ergodic and has
invariant probability density function
\begin{equation}
\label{eq:IM}
\rho_\infty(x) \propto \exp\left(-\frac{1}{\sigma}\Xi(x)\right),
\end{equation}
with respect to Lebesgue measure.

\subsection{Conditional Expectation and Variance}
\label{ssec:CE63}

We aim to predict the $x$ variable from historical data of a long 
trajectory of $x$ alone. Thus the observation and observable maps are
$\ix(\omega) = x, F(\omega) = x.$
We will also estimate second moments, enabling us to compute 
conditional variance, for which $F(\omega) = x^2.$
Observation data is generated by using an implicit time-stepping scheme with
time-step $0.01$ in the slow variable and built-in Matlab solvers to integrate the 
fast variables with $\lorenzeps=0.001$. Source data for the slow variable 
$x$ is gathered for $N=40000$ points sampled at the macroscopic time interval 
$\Delta t = 0.05$. Then $\hat{N}=7500$ out-of-sample points from a new 
trajectory $\{\xo_{n}\}$ are gathered at the same resolution, 
which provide the ground 
truth for assessing forecast error. The natural error metric is the root 
mean squared error (\RMSE), the $L^2$ norm of $Z_\tau - U^\tau F.$ To account for differences in scale, we normalize the \RMSE by the standard deviation of the trajectory.
\begin{displaymath}
\RMSE(\tau) = \bigg(\frac{1}{\hat{N}-q}\sum_{n=0}^{\hat{N}-q-1} |Z_\tau(\xo_n) - \xo_{n+q}|^2\bigg)^{1/2}\bigg(\frac{1}{\hat{N}-q}\sum_{n=0}^{\hat{N}-q-1}|\hat{x}_{n+q}-\bar{x}|^2\bigg)^{-1/2}.
\end{displaymath}
\begin{figure}[t!]
\centering
\begin{overpic}[width=.6\textwidth]{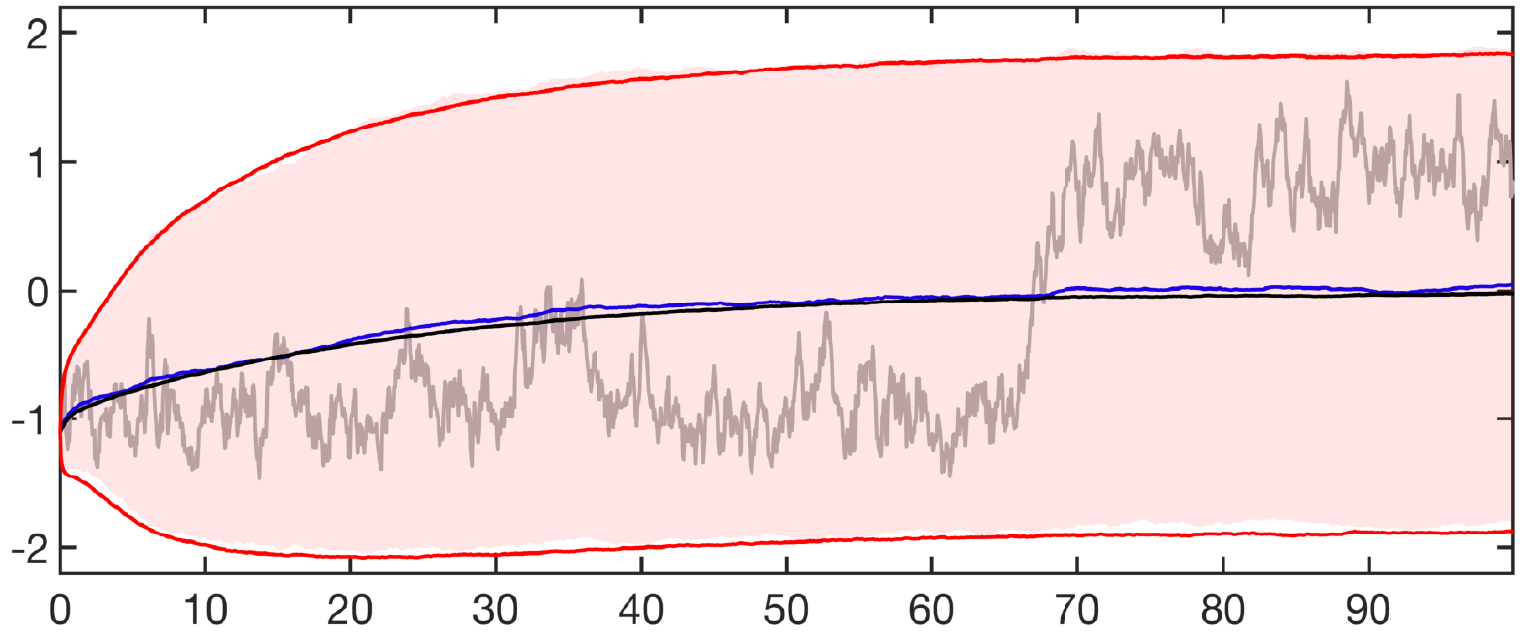} \put(45,-1) {$\tau$}
\put(100,30){\fbox{%
  \parbox{1.9cm}{\scriptsize
    {\color{gray} \bf --}  Trajectory \\
    {\color{black} \bf --}  SDE mean \\
    {\color{red} \bf --}  SDE 2$\sigma$ \\
    {\color{blue} \bf --}  Prediction\\
    {\color{pink}$\blacksquare$} Predicted 2$\sigma$ 
  }%
}}
\end{overpic}
\caption{{\bf Long-term forecast convergence}. Grey: the trajectory of the SDE
started from $x=-1.10$; blue the KAF predictor $Z_\tau(-1.10)$ with pink 
shades giving two standard deviation confidence bands computed
from the conditional variance;  black the Monte Carlo approximation of the 
conditional mean using the SDE; red the Monte Carlo approximation of the 
conditional variances using the SDE. KAF computations of mean and variance agree with the true conditional expectation 
and mean computed from $10000$ Monte Carlo realizations of the SDE.} 
\label{fig:forecastL63}
\end{figure}
\begin{figure}[!ht]
\centering
\begin{overpic}[width=.6\textwidth]{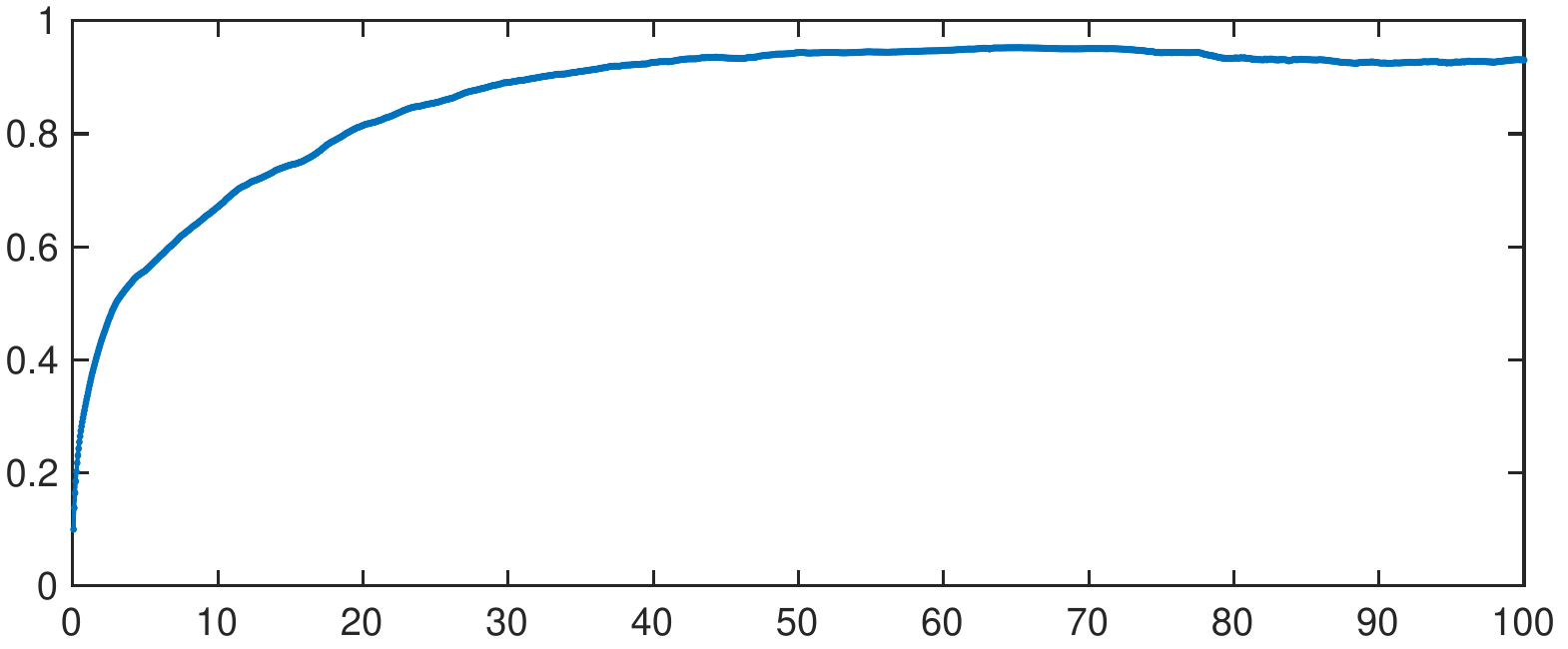} \put(45,-1) {$\tau$} 
\end{overpic}
\caption{{\bf Long-term \RMSE} for the forecast in Figure~\ref{fig:forecastL63} saturates at $\tau=50$, as the forecast converges to the long-term mean at $X=0$. }
\label{fig:rmseL63}
\end{figure}
\begin{figure}[!ht]
\centering%
\begin{overpic}[width=.6\textwidth]{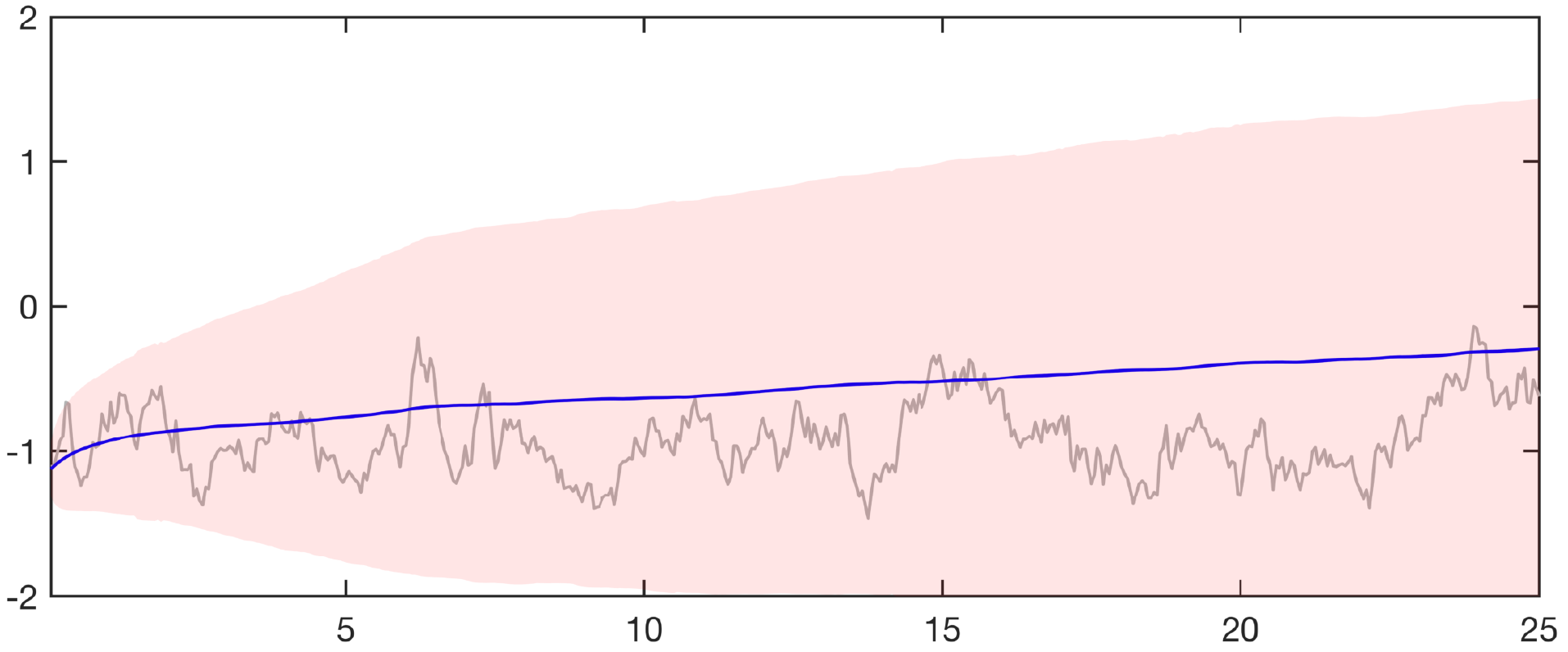} \put(100,33){\fbox{%
  \parbox{1.9cm}{\scriptsize
    {\color{gray} \bf --}  Trajectory \\
    {\color{blue} \bf --}  Prediction\\
    {\color{pink}$\blacksquare$} Predicted 2$\sigma$ 
  }%
}}
\end{overpic}
\begin{overpic}[width=.6\textwidth]{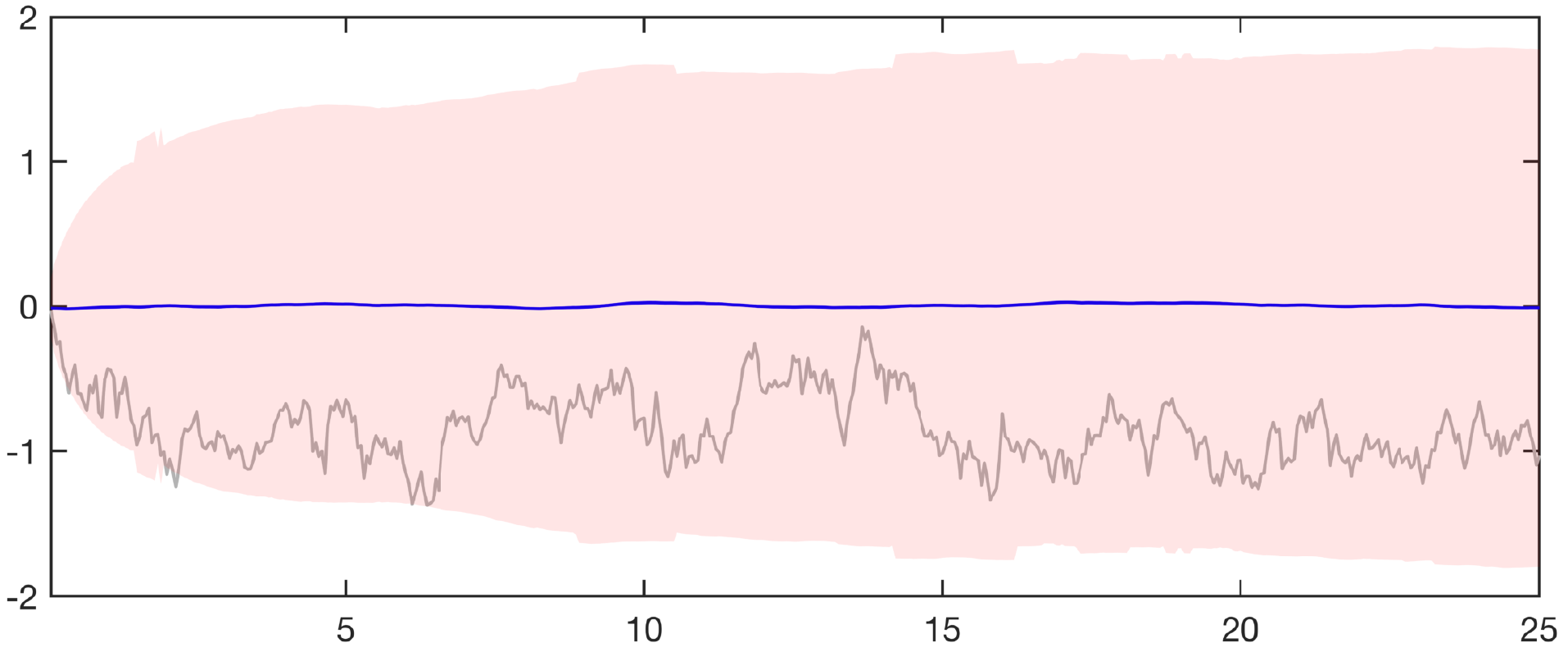} \put(50,0) {$\tau$} \end{overpic}
\caption{{\bf Comparison of uncertainty for different initial data}. 
Blue is predictor (conditional mean), 
grey is trajectory and pink gives two standard
deviations bounds computed from the conditional variance. The forecast 
uncertainty when initialized at $x=0$ (right) displays more rapid
growth in uncertainty over short lead times than that when $x=-1.10$ (left). 
This sensitivity to initial conditions is a desirable feature of KAF, in light of the fact that both plots share the same training set.}
\label{fig:l63_uq_comp}
\end{figure}

Figure~\ref{fig:forecastL63} depicts the behavior of the KAF 
forecast $Z_\tau(x)$ as a function of lead time $\tau$, for
fixed $x$. This forecast exhibits
two interesting properties which can be understood through 
the small-$\lorenzeps$ homogenization limit. The first relates
to the fact that the trajectory itself is not well-predicted;
the second explains what is well-predicted.

Firstly, the predictor 
tracks the conditional mean initialized at $x=-1.10$, and not the
trajectory itself. This is predicted by the theory, since what is
predicted is the long-term conditional expectation 
$\Expect[U^\tau F| \ix]$. Indeed this latter quantity necessarily 
converges for large $\tau$ to a constant, under mixing assumptions 
on the $(x,y)$ system, whilst individual trajectories in $x$ exhibit 
stochastic dynamics, approximately that of $X$. 
This explains the growth in \RMSE seen in Figure~\ref{fig:rmseL63}.
Secondly, exploiting the fact that we expect the system \eqref{eq:l63}
to behave like \eqref{eq:SDE}, when projected onto the $x$ co-ordinate, 
we can provide an objective evaluation of the KAF forecasts by running 
Monte-Carlo simulations of the SDE for $X$; to do this we 
compute the sample mean and variance over $10000$ sample paths 
initialized at the same initial point $x= -1.10$. 
Figure~\ref{fig:forecastL63} compares the KAF forecast mean and variance
with that predicted by  Monte-Carlo mean and variance for the SDE, and they
are seen to agree very well over the entire window of computation. 

Finally in Figure~\ref{fig:l63_uq_comp} we use the possibility of 
varying the initial condition in the KAF to demonstrate that 
the variability encapsulated in the conditional variance is 
able to pick-up different sensitivities, depending upon initial condition. 
The panel on the left shows a trajectory initialized at $x=-1.10$ and the
panel on the left shows a trajectory initialized at $x=0.$
As can be expected from the limiting SDE, the uncertainty when
initialized at $x=0$ is greater and this is manifest in the conditional
variance.

Examination of the KAF technique in this homogenization setting
thus clearly reveals the inability of the method to predict trajectories,
but shows that it can accurately approximate statistics of trajectories,
averaged over the unobserved component of the system.
Furthermore, analysis of the SDE provides a means of characterizing 
the geometry of the underlying hypothesis space, as seen in the next
subsection.

\subsection{Insights Into The Hypothesis Space}
\label{ssec:E}

\begin{figure}[!t]
\centering
\begin{overpic}[width=.31\textwidth]{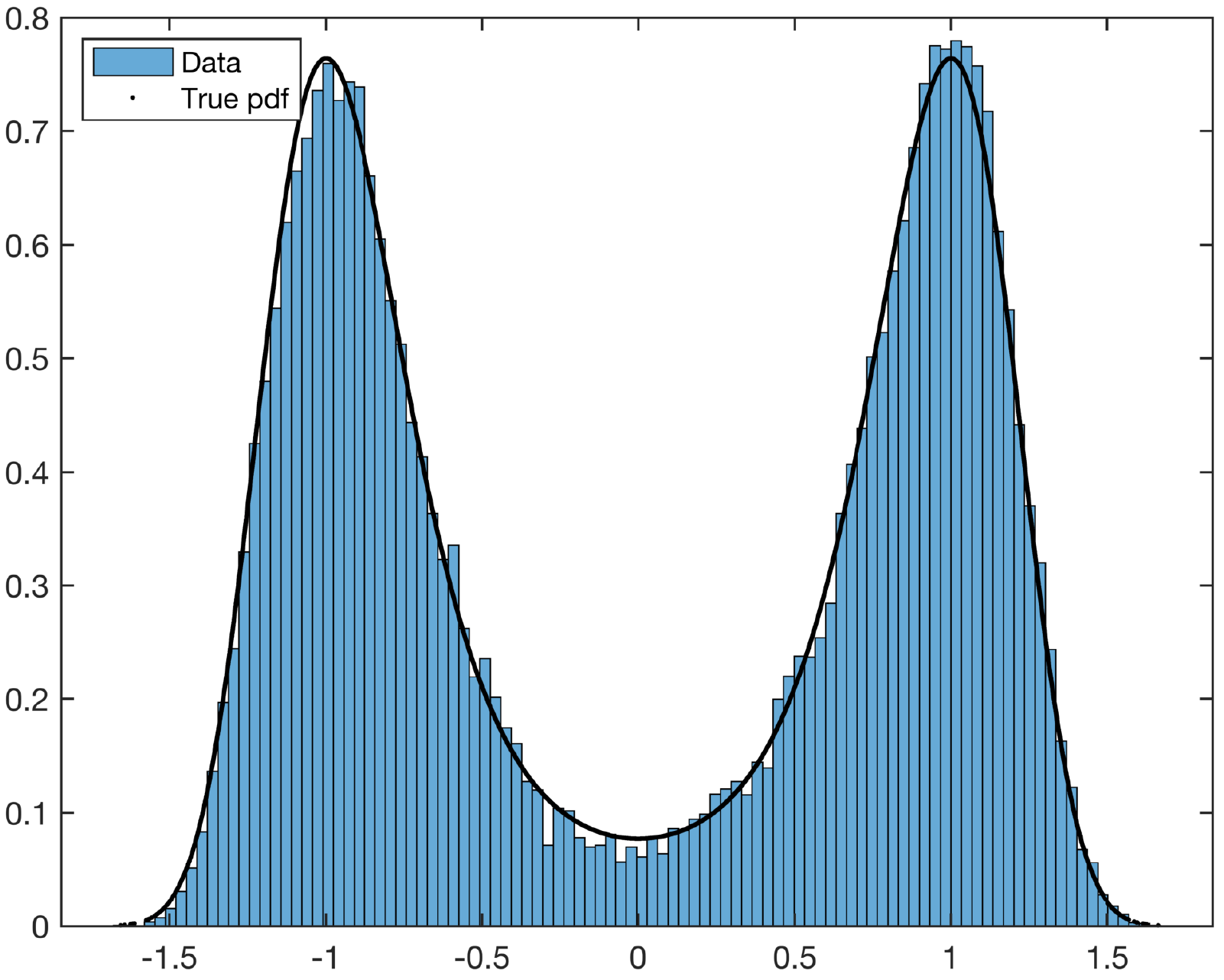} \put(40,80) {\scriptsize $\rho_\infty(x)$} 
\end{overpic}~
\begin{overpic}[width=.32\textwidth]{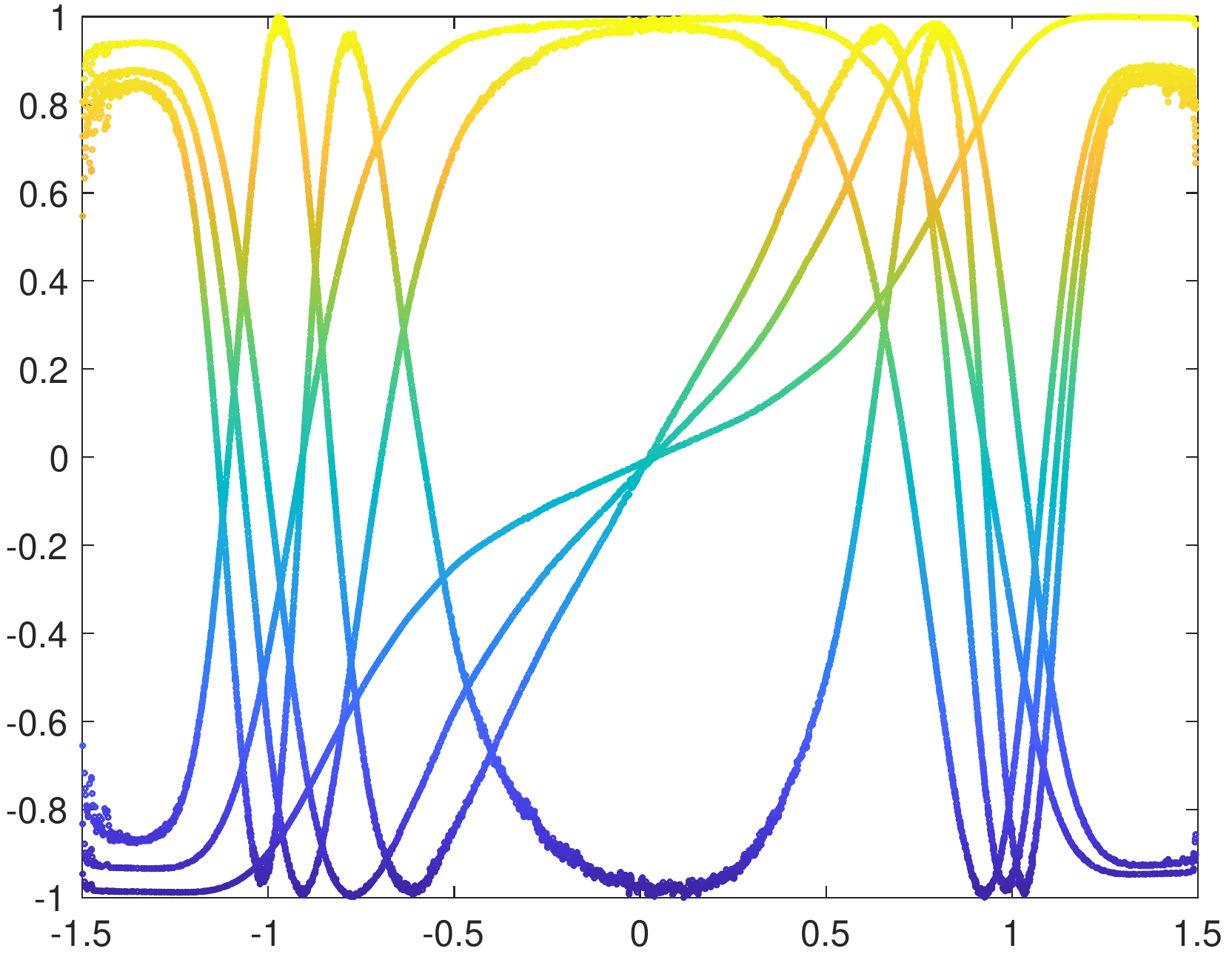} 
\put(40,80) {\scriptsize $\psi_j(x)$}
\end{overpic}~
\begin{overpic}[width=.32\textwidth]{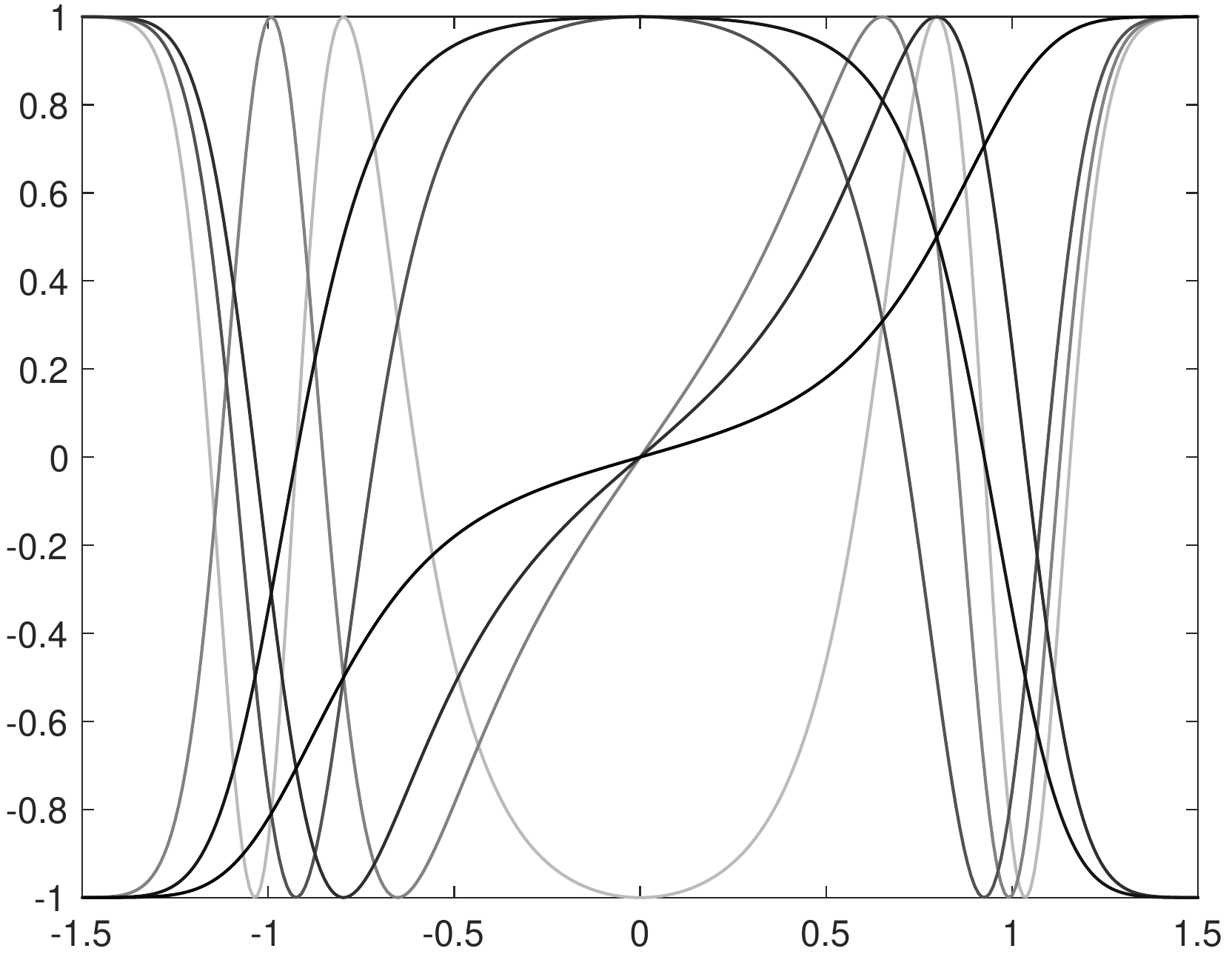} 
\put(40,80) {\scriptsize $\varphi_j(X)$}
\end{overpic}
\caption{{\bf Comparison of invariant densities and eigenfunctions}.
Left: invariant density \eqref{eq:IM}
and histogram of $x$ from \eqref{eq:l63}. Middle and right: six empirically
computed eigenfunctions using $x$ from \eqref{eq:l63}, and 
using theory associated with \eqref{eq:SDE}, respectively.}
\label{fig:phi_compare}
\end{figure}

By studying the large data and small kernel bandwidth limit of the 
matrix $G$ defined in \eqref{eq:G} we get insights into
the structure of the hypothesis space \eqref{eq:HS}. 
The theory in \cite{Giannakis19}, building on the
papers \cite{coifman2006diffusion,BerryHarlim16}, demonstrates that,
for the choice of kernel described in Appendix~\ref{appComputation}, the vectors of $G$ are approximated in the large data and 
small bandwidth limit by eigenfunctions of  
the Laplace-Beltrami operator $\Delta_h$ on $M$, associated 
to a metric $h$. For the diffusion 
process \eqref{eq:SDE} in dimension $m=1$ the manifold $M$ is simply $\R$
and the  metric $h$ is given by $h=\rho^{-2}$, with invariant density $\rho$
given by \eqref{eq:SDE}. 
The fact that the density $\rho$ is approximately available to
us through the time series $x$ generated by \eqref{eq:l63} 
is demonstrated in the left panel of
Figure \ref{fig:phi_compare}, where we compare the histogram
generated by the data with $\rho$ given by \eqref{eq:IM}.
The conclusion of these various approximations is that we expect
the eigenvectors of $G$, based on data $x$ from \eqref{eq:l63}, to
be well-approximated by eigenfunctions of $\Delta_h$ on $\R$, with
$h=\rho^{-2}.$ We now demonstrate that this is indeed the case.

The action of the Laplace-Beltrami operator on $f\in C^\infty(M)$ is 
given by
\begin{displaymath}
    \Delta_h f = -{\divr}_{\mu}{\grad}_h f = -{\divr}_{\mu}\left(\frac{1}{\rho} \frac{df}{dX}\right) =  -\frac{1}{\rho}\frac{d}{dX}\left(\frac{1}{\rho}\frac{df}{dX}\right), 
\end{displaymath}
where $\divr_\mu$ and $\grad_h$ are the divergence and gradient operators associated with $\mu$ and $h$, respectively. Using the above, we solve the eigenvalue problem $
    -\Delta_h \varphi = \lambda^2\varphi$ directly. 
We make the substitution $dY=\rho dX$, mapping $X \in \R$ into $Y \in [0,1]$;
we note that $Y$ has interpretation as the cumulative distribution function 
coordinate of $\rho$. In terms of $Y$ we have
\begin{displaymath}
    \Delta_h f =  -\frac{1}{\rho}\frac{d}{dX}\left(\frac{1}{\rho}\frac{df}{dX}\right) = -\frac{d^2f}{dY^2} .
\end{displaymath}
Noting that the natural boundary conditions for the Laplace-Beltrami
operator are of no-flux type, it follows that, when viewed as functions of
$Y$, the eigenfunctions of $\Delta_h$ 
satisfy a boundary value problem of the form 
\begin{align*}
    -\varphi''(Y) = \lambda^2\varphi(Y), \\
    \varphi'(0) = \varphi'(1) = 0. 
\end{align*}
The solutions are the well-known harmonics
$\cos(k\pi Y),$ and corresponding eigenvalues 
$\lambda_k = k\pi, k\in\mathbb{N}$. 
Changing back to variable $x$ we obtain
\begin{align}
\label{eq:varphi}
    \varphi_k(X) = \cos \left(k \pi \int_{-\infty}^X \rho(z) dz\right).
\end{align}
We now verify that, for large data sets and small
bandwidth, the eigenfunctions of $G$ are indeed close to those
associated with the Laplace-Beltrami operator $\Delta_h.$ 
This is demonstrated in the middle and right panels of
Figure~\ref{fig:phi_compare}. The middle panel shows the first
six eigenfunctions of $G$, computed from data derived from the
$x$ variable in \eqref{eq:l63}; the right panel shows the
eigenfunctions of $\Delta_h$ for diffusion process \eqref{eq:SDE} in
variable $X$ which governs the limiting behaviour of $x$ in the 
scale-separated case. The agreement is excellent, demonstrating
that the heuristics underlying parameter choices within the
kernel based methodology (see Appendix~\ref{appComputation})
work well in this setting.

\begin{figure}[!ht]
\centering
\begin{overpic}[width=.6\textwidth]{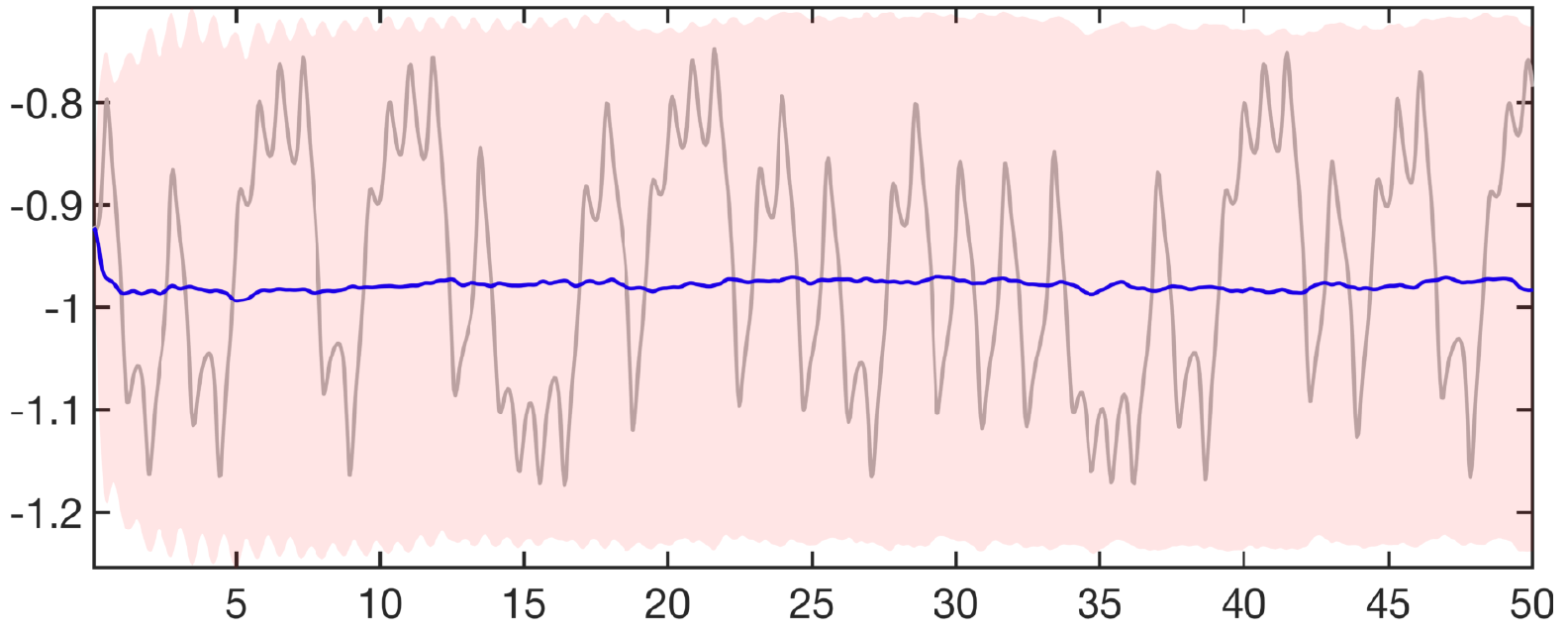} 
\put(40,37) {\footnotesize $\ix(\omega) = x$} 
\put(100,33){\fbox{%
  \parbox{1.9cm}{\scriptsize
    {\color{gray} \bf --}  Trajectory \\
    {\color{blue} \bf --}  Prediction\\
    {\color{pink}$\blacksquare$} Predicted 2$\sigma$ 
  }%
}}
\end{overpic}

\begin{overpic}[width=.6\textwidth]{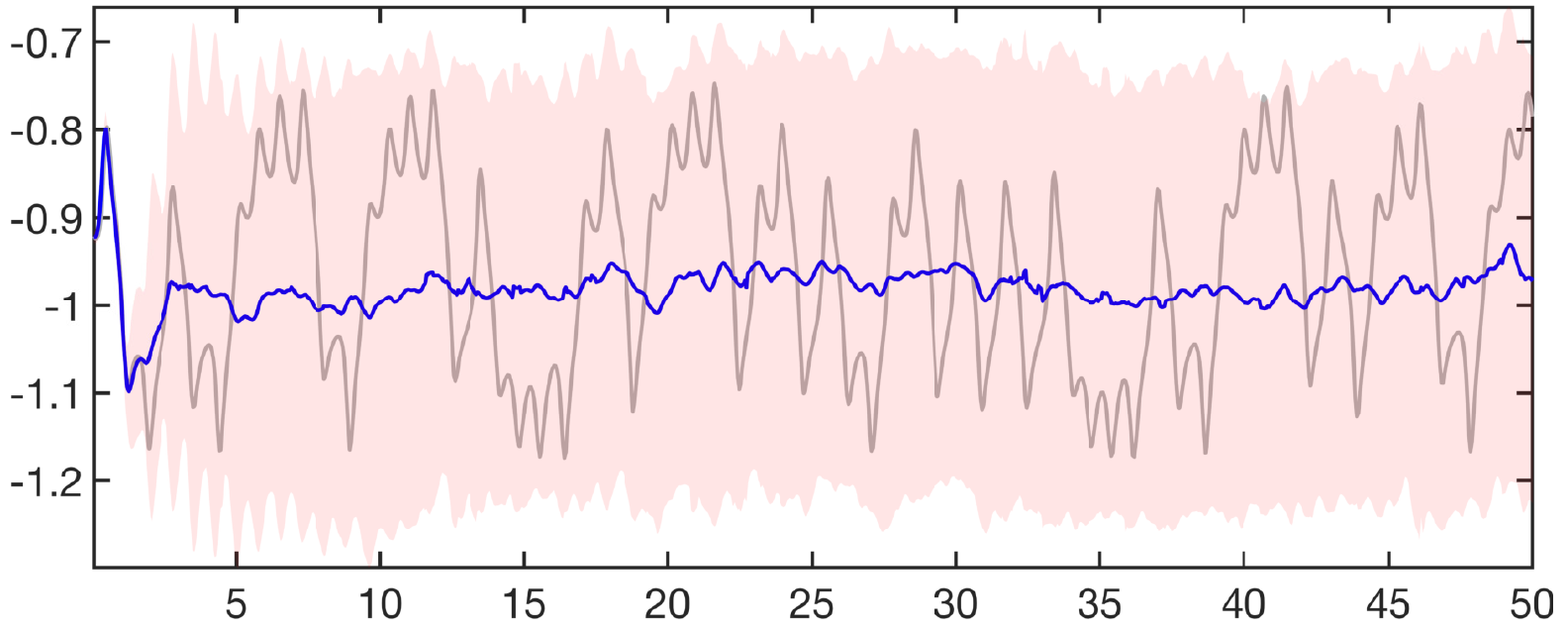} 
\put(50,-2) {$\tau$}
\put(32,35) {\footnotesize $\ix(\omega) = [x,4y_2/90]^T$} 
\end{overpic}
\caption{{\bf Prediction in non-Markovian regime, $\lorenzeps= 1$, } requires observations augmented by the forcing term to achieve short-term predictability.}
\label{fig:nonMarkovL63}
\end{figure}

\subsection{Non-Markovian Regime}
\label{ssec:NM1}

In the preceding subsections we studied predictors for $x$, based only
on time-series data in the $x$ coordinate, for the equation \eqref{eq:l63}.
We studied the scale-separated regime where $\lorenzeps \ll 1$ and $x$
is approximately Markovian -- it is approximately governed by an
SDE. Here we study the behavior of identical 
predictors when $\lorenzeps=1$; the system ~\eqref{eq:l63} then no longer 
exhibits homogenization and $x$ is no longer Markovian in view of
the lack of scale-separation. As a result, 
the prediction of $x$, shown in the top of Figure~\ref{fig:nonMarkovL63}, 
is poor even at very short times, and the two standard deviation
confidence bands reflect this rapid initial error growth, 
and then remain large throughout the time window. 
$Z_\tau$ converges rapidly to the conditional mean. 
To render the problem Markovian we may include more data, and
in particular the forcing term. To this end we take 
$\ix(\omega) = [x,4y_2/90]^T$; see the bottom of 
Figure~\ref{fig:nonMarkovL63}. The prediction of $x$ with these 
augmented observations yields very accurate predictions over 
a lead time of $\tau=3$, considerably larger than in the preceding
case where $\ix(\omega) = x$.  After $\tau=3$, however,
pathwise predictability again fails, due to the sensitive dependence of
solutions with respect to the forcing function. Once again
$Z_\tau$ converges to the conditional mean.
However the uncertainty quantification provided by two
standard deviation error bars 
consistently captures the true trajectory, even in this 
non-Markovian regime.

This non scale-separated pair of examples illustrates that 
prediction of non-Markovian variables is inherently harder than
Markovian variables, but that sensitive dependence of trajectories
with respect to perturbations also limits predictability,
even in the Markovian setting. This second point will be
prominent, too, in the next section.

\section{Averaging: Lorenz 96 Multiscale System}

This section is devoted to the setting in which a chaotic ODE
of form (A) is approximated by an ODE of form (A0).
The goal is to make predictions of the $x$ variable using data
concerning the $x$ variable alone from (A). 
The role of (A0) is primarily to help us interpret those
predictions; however it also serves to motivate a different 
prediction methodology,
which mixes model-based and data-driven methodologies, against
which we will compare KAF. 

The averaging setting presents interesting 
opportunities to understand forecasting. In particular, by tuning a
parameter in (A), which is also present in (A0), we are able
to create settings in which the variable to be predicted is,
in turn, approximately periodic, quasi-periodic and chaotic. 
These different settings give rise to different types of forecasts 
and we demonstrate this. In subsection~\ref{ssec:CE96} we examine 
the long-term statistics predicted by KAF from data generated 
by (A) -- the conditional expectation
and variance of one component of the slow variable --
and compare them with estimates
computed from (A0) using Monte-Carlo simulation. This illustrates
our main contribution 1 from the list in subsection \ref{ssec:OC}.
Then, in subsection~\ref{ssec:GP96}, we compare the purely
data-driven method of prediction with a hybrid data-model
predictor. This hybrid is computed by using  Gaussian process
(GP) regression to compute an approximate closure $\hgp(x) \approx v(x)$, 
in the notation of (A0); for more details on such approximate
closures in the context of the specific model we study in
the following four subsections
see \cite{fatkullin2004computational} and references therein.  
The  work in subsection~\ref{ssec:GP96}  
illustrates our main contribution 3  from subsection~\ref{ssec:OC}.
Subsection~\ref{ssec:NM2} is concerned with non-Markovian prediction,
in which there is no scale-separation between observed and unobserved
variables.
We start, however, in subsection~\ref{ssec:L96}, introducing the
concrete model around which our experiments are organized.

\label{sec:96}

\subsection{The Model}
\label{ssec:L96}

\begin{figure}[t]
\centering
\begin{overpic}[width=.3\textwidth] {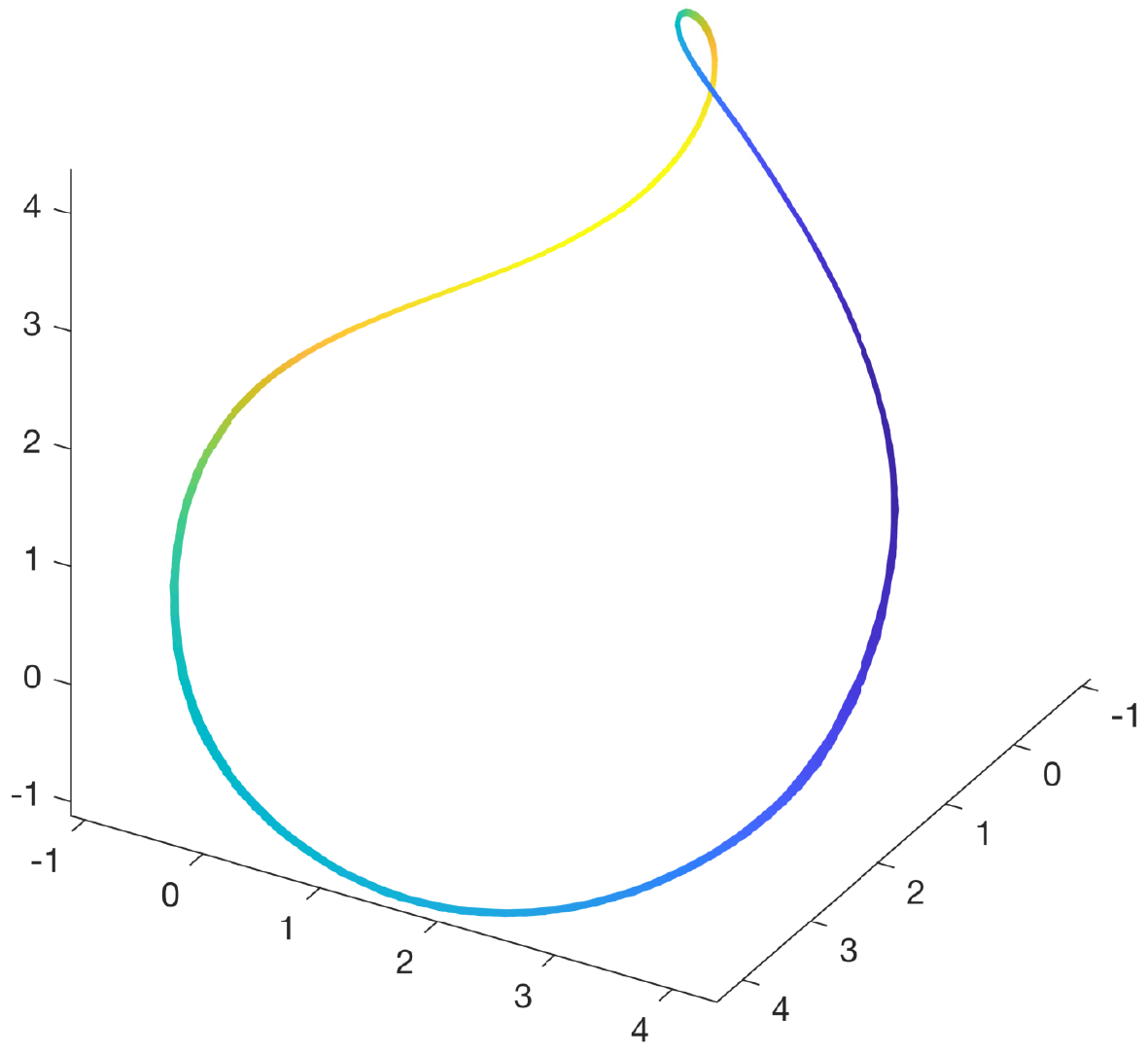} 
\put(5,80) {\footnotesize Periodic} 
\end{overpic}
\begin{overpic}[width=.3\textwidth]{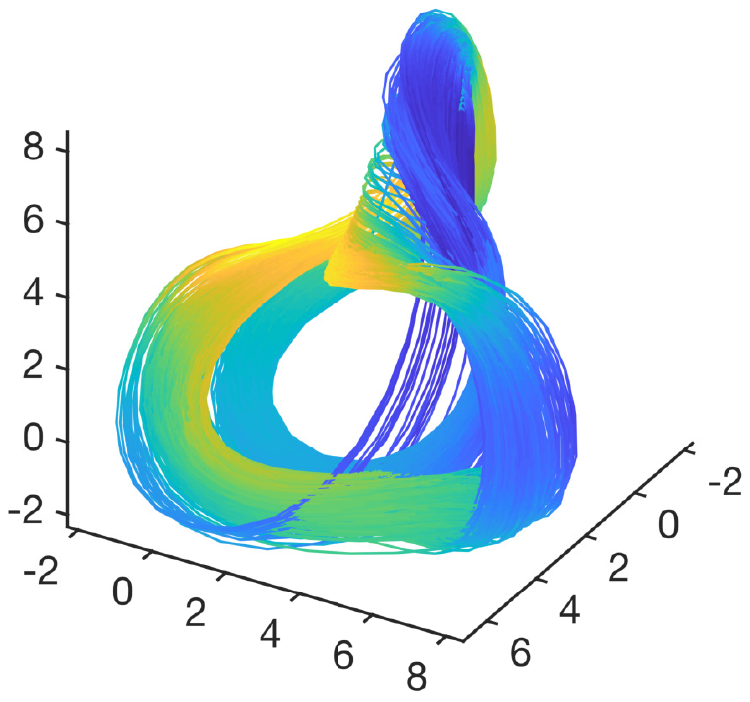} 
\put(5,80) {\footnotesize Quasiperiodic} 
\end{overpic}\begin{overpic}[width=.3\textwidth]{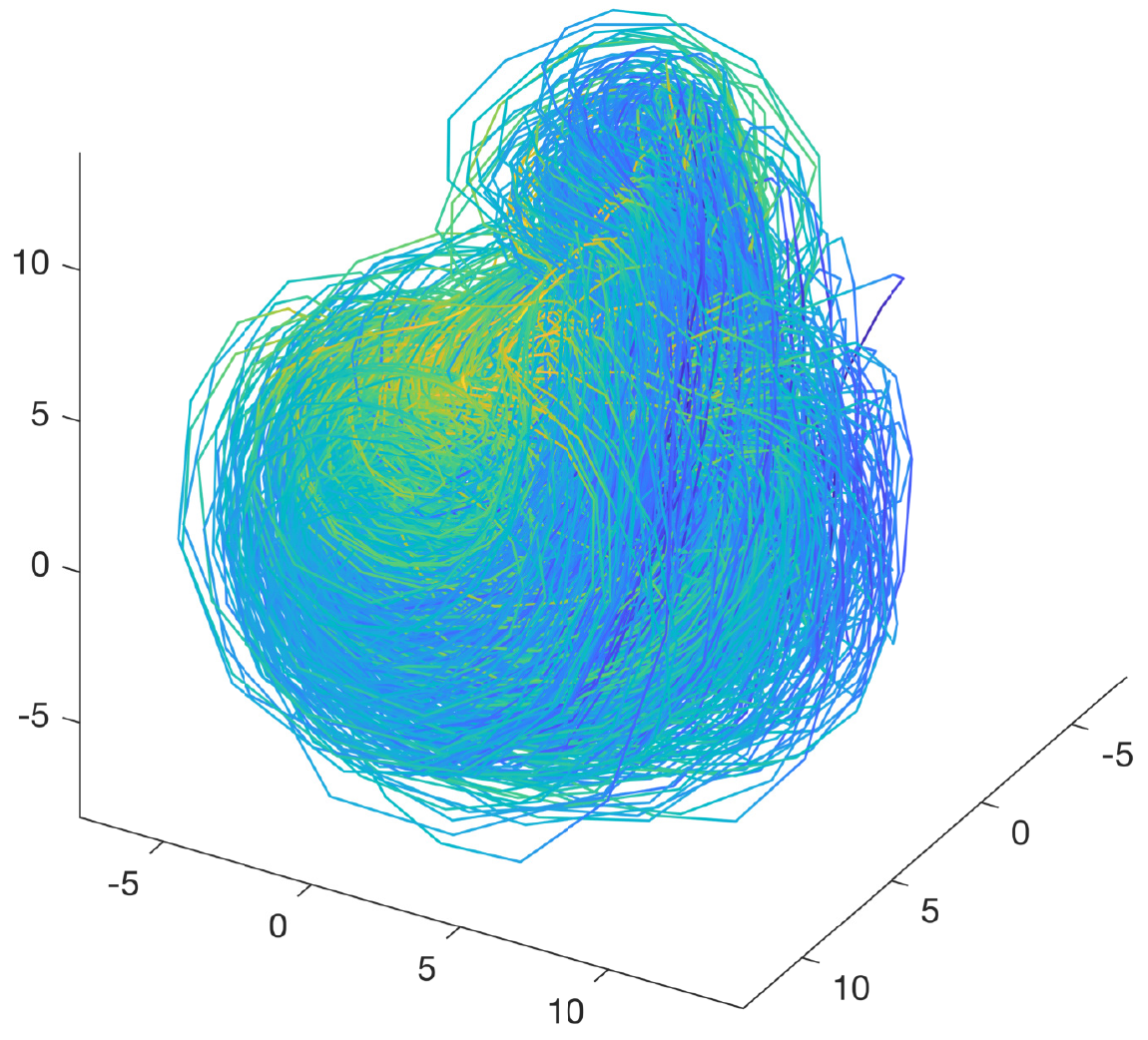} 
\put(5,80) {\footnotesize Chaotic} 
\end{overpic}
\caption{
  \textbf{Lorenz 96 regimes of increasing complexity (left to right)}.
  Phase portraits show $(x_1,x_2,x_3)$ coordinates shaded by $x_4$.
The parameter $F_x$ takes values $5.0, 6.9$ and $10.0$ respectively,
from left to right, and all other parameters are as in \eqref{eq:param}.}
\label{fig:phase_l96}
\end{figure}

In this section we focus our attention on another chaotic dynamical system,
colloquially known as ``Lorenz 96 multiscale''~\cite{lorenz1996predictability},
which we will simply abbreviate to \LNS{}. Following the notation established
in~\cite{fatkullin2004computational}, the \LNS{} equations model $K$ slow variables $\{x_k\}_{k=1}^K$ coupled to $JK$ 
fast variables $\{y_{j,k}\}_{j,k=1,1}^{J,K}$ 
with evolution given as follows:
\begin{equation}
\label{eq:l96}
  \begin{aligned}
    \dot{x}_k&=-x_{k-1}(x_{k-2}-x_{k+1})-x_{k}+F_x+\frac{h_x}{J}\sum_{j=1}^Jy_{j,k},\\
    \dot{y}_{j,k}&=\frac{1}{\lorenzeps}\Bigl(-y_{j+1,k}(y_{j+2,k}-y_{j-1,k})-y_{j,k}+h_y x_k\Bigr),\\
&\quad x_{k+K}   = x_k, \quad y_{j,k+K} = y_{j,k}, \quad y_{j+J,k} = y_{j,k+1}.
  \end{aligned}
\end{equation}
This is of the form (A).
On the assumption that the $y$-variables, with $x$ frozen, are ergodic, the
averaging principle shows the existence of a function $C\colon \R^K \to
\R^K$ such that, for small $\lorenzeps$, the $x$ variables are approximated by
$X=(X_1, \dots, X_k)$ solving
\begin{equation}
  \label{eq:l96a}
  \dot{X}_k = -X_{k-1} (X_{k-2} - X_{k+1}) - X_{k} + F_x + h_x C_k(X), \quad k \in \{1,\dots, K\},
\end{equation}
with the periodic boundary conditions $X_{k+K} = X_k$ and
$C_k\colon \R^K \to \R$ denoting the $k^{th}$ component of vector-valued
function $C$.  This system is of form (A0).
Since system~\eqref{eq:l96} is index-shift-invariant, it is clear that the
closure $C_k$, if it exists, satisfies $C_{k+1}(X)=C_k(\pi X)$ where $\pi$
shifts the vector indices by adding one unit, 
invoking periodicity at the end points. 
Furthermore, when $J$ is large, \revisionOne{empirical evidence
\cite{wilks2005effects,fatkullin2004computational}} suggests that there is a
function $c\colon \R \to \R$ such the approximation $C_k(X)=c(X_k)$ is a good
one. 
For the current work the key point to appreciate is that for small $\lorenzeps$
the variables $x$ in \eqref{eq:l96} exhibit (approximately) Markovian behavior,
and this behavior is deterministic and governed by $X$. However, by
tuning $F_x$, different responses arise in the deterministic variable.
In the following  we fix parameters $\lorenzeps^{-1},K,J, h_x, h_y$ 
throughout all our experiments as follows:
\begin{equation}
\label{eq:param}
\lorenzeps^{-1} = 128, \, K = 9, \, J = 8, \, h_x = -0.8, \, h_y = 1.0.
\end{equation}
We then choose $F_x$ to distinguish
three cases as follows:
\[
  \begin{gathered}
    \text{periodic} \\
    \begin{aligned}
      F_x &= 5.0,
    \end{aligned}
  \end{gathered}
  \qquad \quad
  \begin{gathered}
    \text{quasi-periodic} \\
    \begin{aligned}
      F_x &= 6.9,
    \end{aligned}
  \end{gathered}
  \qquad \quad
  \begin{gathered}
    \text{chaotic} \\
    \begin{aligned}
      F_x &= 10.0.
    \end{aligned}
  \end{gathered}
\]
Figure \ref{fig:phase_l96} demonstrates the three responses within
system \eqref{eq:l96} resulting from these parameter choices.

\subsection{Conditional Expectation and Variance}
\label{ssec:CE96}

We aim to predict the $x_1$ variable from historical data of a long trajectory of $x$ alone. Thus the observation and observable maps are
$\ix(\omega) = x, \, F(\omega) = x_1.$
We will also use $F(\omega) = x_1^2$ when estimating conditional variance.
By tuning the scalar parameter $F_x$ (not to be confused with function $F$)
as outlined in the preceding subsection we can obtain periodic, quasi-periodic 
and chaotic responses in the averaged variable $X$.
It is intuitive that the ability of the KAF to track the true trajectory of
the slow variables decreases with increasing complexity; in other words,
predictions in the periodic case should be the most
accurate whilst those in the chaotic case present a significant challenge.
In the experiments that follow the size and sampling interval of the source 
(training) data remain fixed at $(40000,0.05)$ and the out-of-sample (test) 
data set is fixed at $\hat{N}=7000$.

\begin{figure}[t]
\centering
\begin{overpic}[width=.6\textwidth]{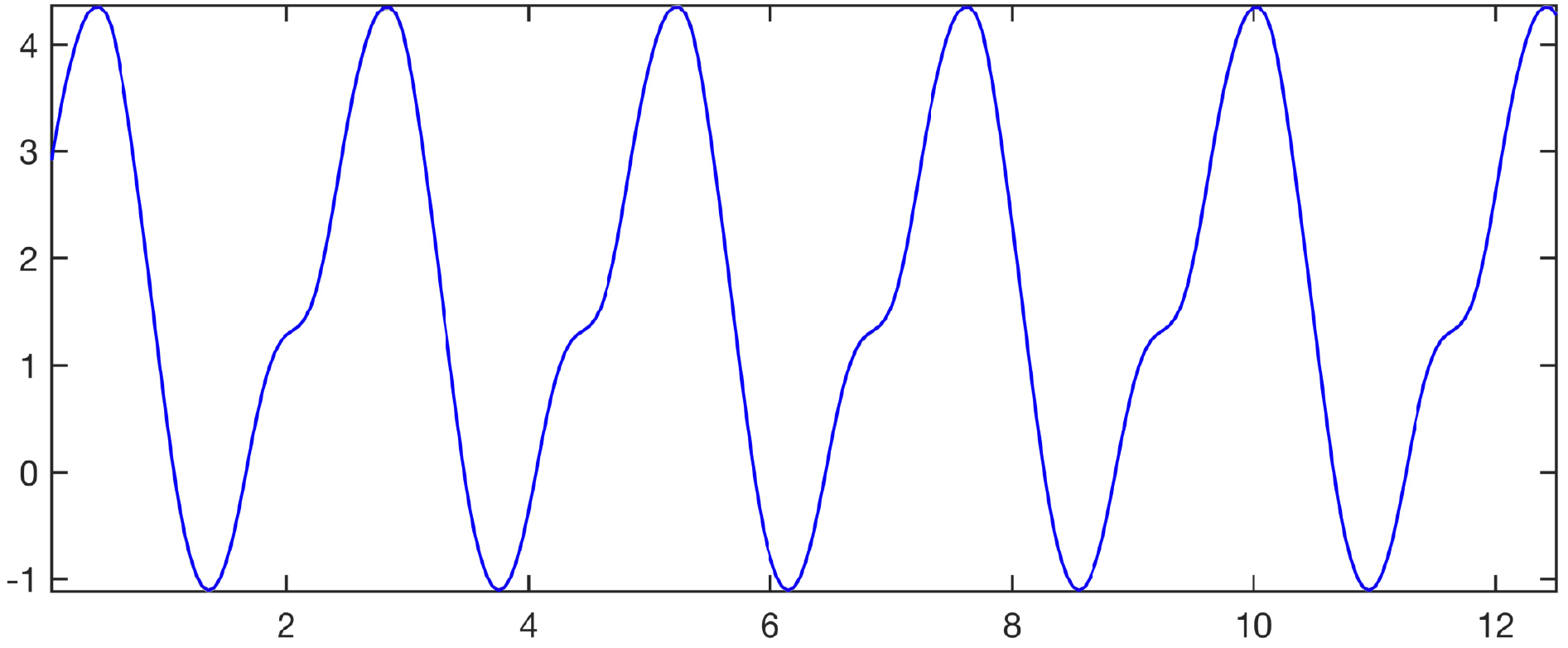} 
\put(8,37) {\footnotesize Periodic} 
\put(100,33){\fbox{%
  \parbox{1.9cm}{\scriptsize
    {\color{gray} \bf --}  Trajectory \\
    {\color{blue} \bf --}  Prediction\\
    {\color{pink}$\blacksquare$} Predicted 2$\sigma$ 
  }%
}}
\end{overpic}

\begin{overpic}[width=.6\textwidth]{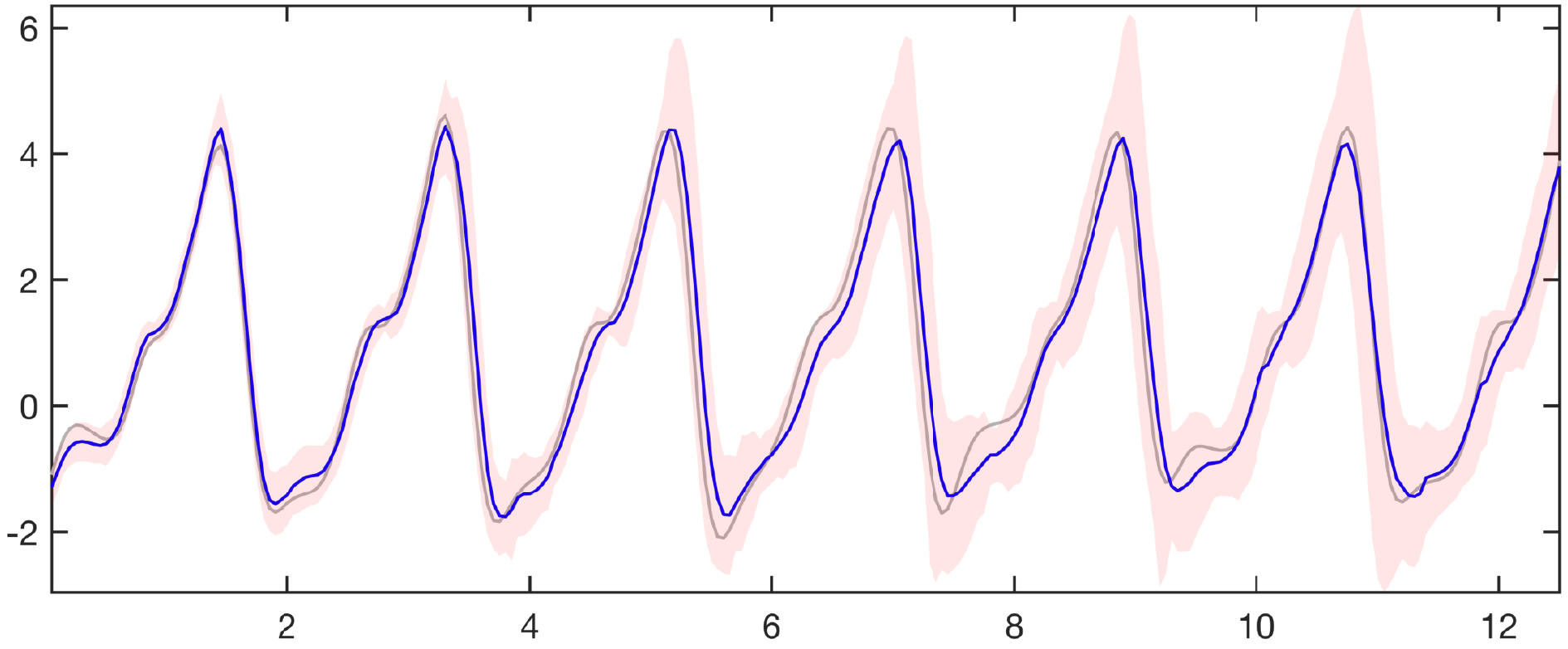} 
\put(8,37) {\footnotesize Quasiperiodic} 
\end{overpic}

\begin{overpic}[width=.6\textwidth]{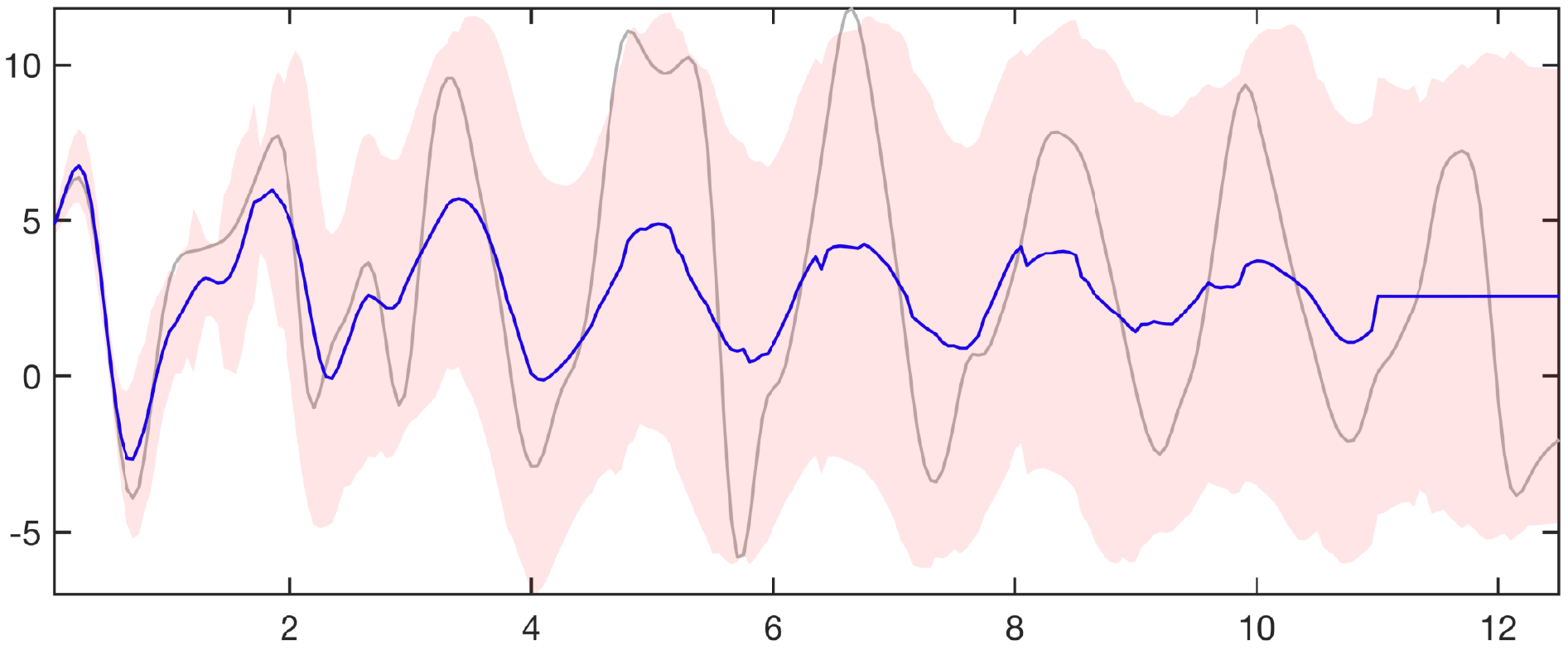} 
\put(50,-2) {$\tau$}
\put(8,37) {\footnotesize Chaotic} 
\end{overpic}
\caption{{\bf Predictability: Periodic, quasi-periodic and chaotic regimes.}
 Prediction $Z_\tau(x)$ of observable $F(\omega) = x_1$ across 3 different regimes.
In each figure grey is the true trajectory, blue the predictor using KAF,
and pink gives two standard deviations confidence bands, computed using the
conditional variance. The parameter $F_x$ takes values $5.0, 6.9$ and $10.0$ 
respectively, from top to bottom, and all other parameters are as 
in \eqref{eq:param}. In the first, periodic response regime, the trajectory is predicted almost
perfectly and this accuracy is reflected in the narrow confidence bands. 
In the second, quasi-periodic response regime, the trajectory 
is predicted very well, but with growing error reflected accurately
in the slowly growing confidence bands. In the third, chaotic response 
regime, the predictive capability is lost due to sensitivity to 
initial conditions and this is reflected in the rapidly growing
confidence bands and in the convergence of the predictor to a constant,
for large $\tau.$}
\label{fig:forecast96}
\end{figure}

The space of observables $\X$ in the current example is the space of all slow
variables.  Since, under the small-$\lorenzeps$ limit, an ODE
closure of the slow dynamics is obtained, 
the variable $x$ behaves (approximately) like a deterministic Markov
process, and the expectation in \eqref{eq:g} disappears; 
the predictor is expected to track
the actual trajectory $x_1(t)$.
To see this another way, note that simply knowing the initial 
values of the $x$-variables (recall that $\X$ is precisely all 
$x$-variables) and the closure $C(X)$
in equation~\eqref{eq:l96a}, we are able to predict 
$x_1$ (or indeed, any $x_k$) exactly, given the initial
conditions for all $x$-variables. 

\begin{figure}[t]
\centering
\begin{overpic}[width=.6\textwidth]{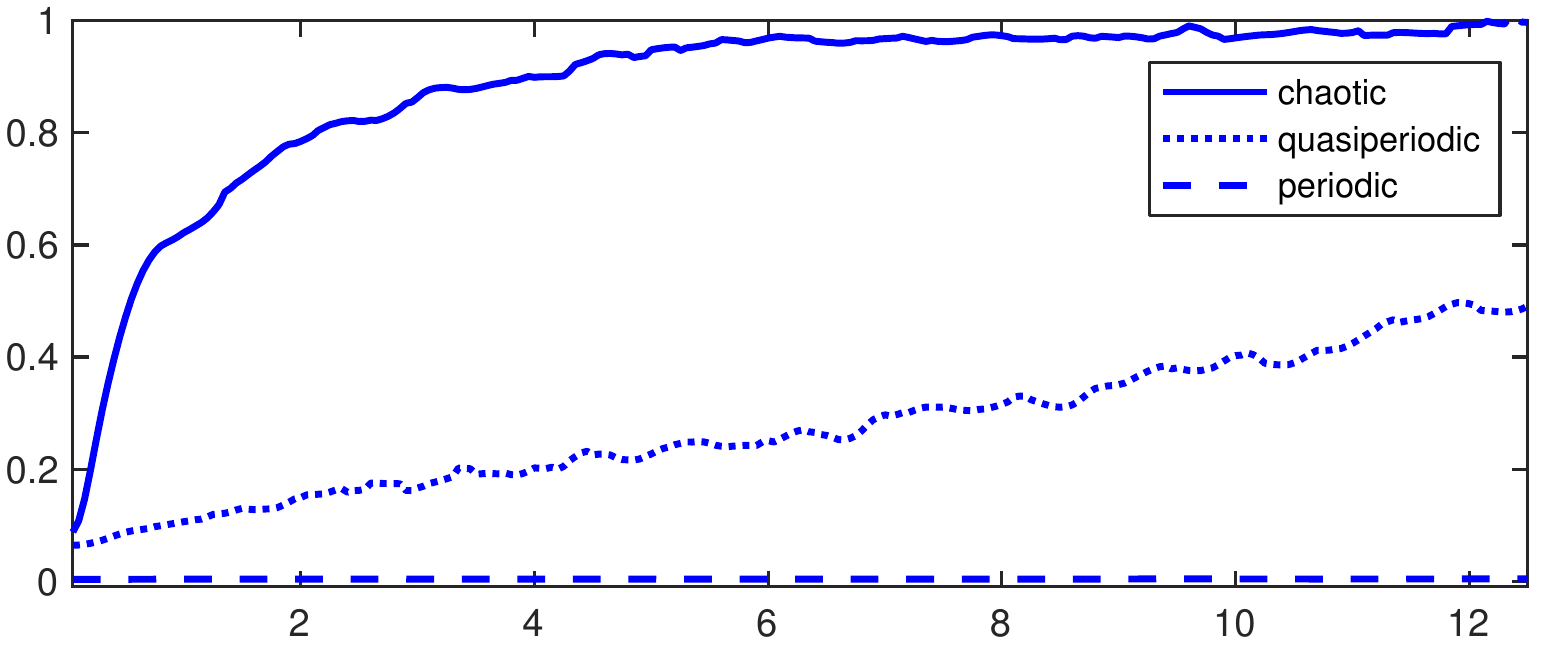} 
\end{overpic}
\caption{{\bf \RMSE}. This figure depicts the RMSE of the predictor $Z_{\tau}(x)$ for \eqref{eq:l96}, for different $F_x$, as a function of $\tau.$
The parameter $F_x$ takes values $5.0, 7.1$ and $10.0$
respectively, from smaller to larger error, corresponding to periodic,
quasi-periodic and chaotic response; all other parameters are as
in \eqref{eq:param}.}
\label{fig:rmse96}
\end{figure}

However this picture is greatly
affected by the sensitivity of the system to initial conditions 
and sampling errors due to high dimensionality of the attractor. 
We now describe how these predictions work in practice, in the
three regimes shown in Figure~\ref{fig:phase_l96}. We display
our results in Figure~\ref{fig:forecast96}, where $x_1$ and 
standard deviation bands are predicted and compared with the true
signal starting from the same point. 
The long-term predictability in each regime is constrained by the complexity of
the underlying Markovian, deterministic, slow dynamics.
In the periodic regime, since chaos is absent in the slow variables, 
a perfect predictor is obtained via the partially observed dynamics;
one interpretation of why this occurs is because \revisionOne{the eigenfunctions of the Koopman operator lie 
in a finite span of the diffusion coordinate observables}~\cite{arbabi2017ergodic}. 
Observe that $Z_\tau$ remains in phase, and the forecast variance is negligible,
for long lead times up to the length of the entire out-of-sample trajectory
($\tau = 350$). 
The quasiperiodic trajectory is tracked imperfectly, but with significant
accuracy over the same range of times; errors are visible
mainly around
the extrema of $x_1$ as suggested by the phase portrait; the conditional
variance reflects the significant accuracy present. 
Prediction in the fully chaotic regime only tracks the trajectory, however,
until a lead time of approximately $1$ time unit, exhibiting behaviour at 
long lead times which is somewhat similar to that seen in the previous, 
homogenization, section in which the predicted variable behaved as if
drawn from a Markov stochastic process. 
In particular the long-term predictor in the chaotic regime converges to a
constant by construction, assuming mixing, and this is consistent with the 
inherent unpredictability of chaotic dynamics.
It is notable that the size of the conditional variance, and the
resulting confidence bands, is a useful guideline as to the
pathwise accuracy of the data-driven predictor.
The observations about the predictability of the system by KAF methods
are also manifest in Figure~\ref{fig:rmse96} which shows the \RMSE in each
of the periodic, quasi-periodic and chaotic regimes.

We mention that in the quasi-periodic case the presence of multiple 
attractors (or multiple lobes of the same attractor), and resulting
intermittent switching between these attractors, leads to a loss of
predictability that is significant on time-scales much longer than
those shown here. For the figure shown here we have ensured that 
training points and out-of-sample points are gathered from the same 
(part of the) attractor to maintain accuracy. We train using two different 
trajectories to gather ample training data.

Recall that at each lead time $\tau$ along the horizontal axis there is
a potentially different number of eigenfunctions $\ell(\tau)$ used in the
data-driven method. See Appendix \ref{as:L} for details on the choice
of $\ell$. In the chaotic regime the optimal $\ell(\tau) $ tends to $1$ 
for large times (see subsubsection \ref{sssec:longt} for an explanation), 
whilst $\ell$ fluctuates around $50$ in the quasiperiodic regime; 
we obtain $\ell\approx 9$ for all $\tau$ in the periodic regime.

\subsection{Comparison Of Data-Driven And Model-Data-Driven Prediction}
\label{ssec:GP96}

The previous subsection concerned purely data-driven prediction of
variable $x$ from (A), using only data in the form of a time-series for $x$.
In this subsection we provide comparison with a different forecasting 
technique based on a combination of model and data-driven prediction,
using data in the form of a time-series for $(x,By)$. Knowledge
of $By$ enables the use of Gaussian process regression (GPR)~\cite{GPML}
(or kriging) to approximate $v(\cdot)$ by $\hgp(\cdot)$ in (A0). Our approach
is motivated by the paper~\cite{fatkullin2004computational} which looked
at finding such closures for the L-96 model in form \eqref{eq:l96}. 
When applied to \eqref{eq:l96} the methodology leads to an approximate
closure for the slow variable $X$ which takes the form
\begin{equation}
  \label{eq:l96b}
  \dot{X}_k = -X_{k-1} (X_{k-2} - X_{k+1}) - X_{k} + F_x + h_x \cgp(X_k), \quad k \in \{1,\dots, K\},
\end{equation}
subject to periodic boundary conditions $X_{k+K} = X_k$. 
This should be compared with \eqref{eq:l96a}, which arises 
from application of the averaging principle; note that, in addition, 
we have invoked the hypothesis that $C_k(X)$ can
be well-approximated by function of $c(X_k)$, as discussed directly 
after~\eqref{eq:l96a}; and we will determine an approximation $\cgp$ for
$c$ by GPR.

\begin{figure}[t]
  \centering
  \includegraphics[width=0.7\textwidth]{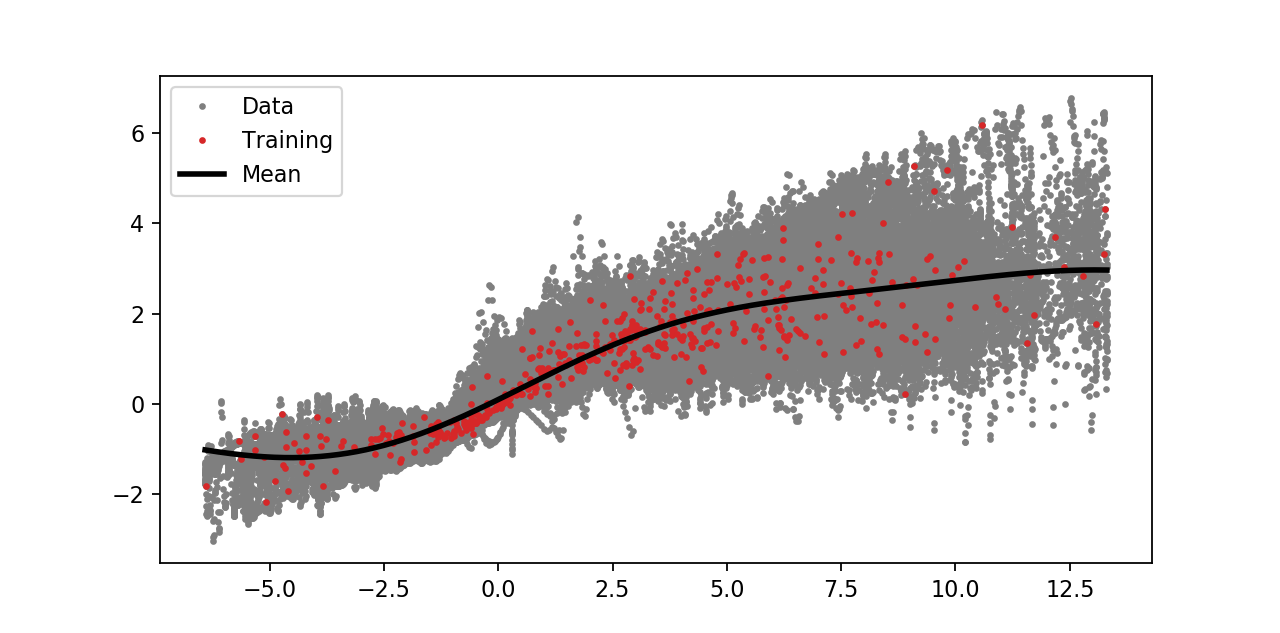}
  \caption{{\bf Mean of Gaussian process regression as a closure.}
Function $\cgp$, and data used to determine it, from data generated
by \eqref{eq:l96} with parameters as in \eqref{eq:param} and $F_x=10.0$}
  \label{fig:gpr}
\end{figure}
\begin{figure}[t]
\centering
\begin{overpic}[width=.6\textwidth]{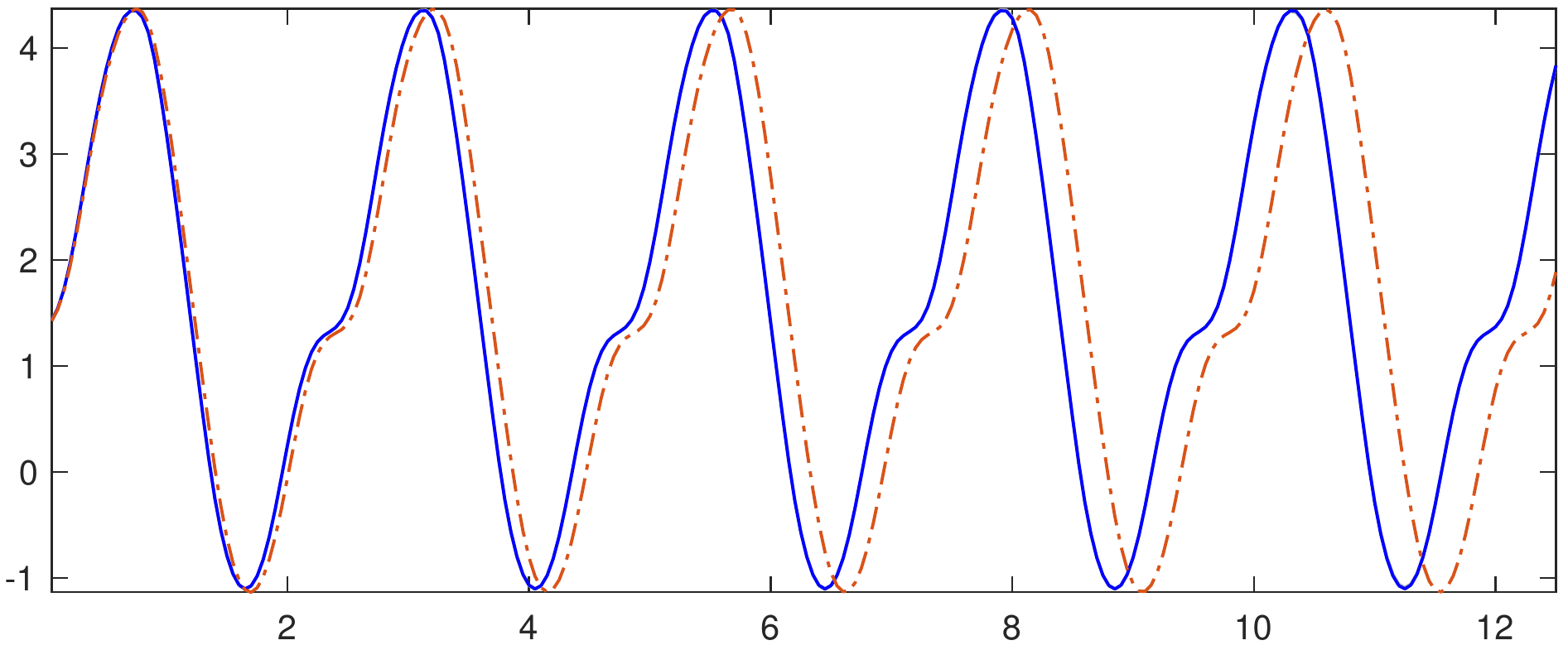} 
\put(8,37) {\footnotesize Periodic} 
\put(100,34){\fbox{%
  \parbox{1.6cm}{\scriptsize
    {\color{gray} \bf --}  Trajectory \\
    {\color{blue} \bf --}  Prediction\\
    {\color{orange} \bf $\cdot$-}  GP closure
  }%
}}
\end{overpic}

\begin{overpic}[width=.6\textwidth]{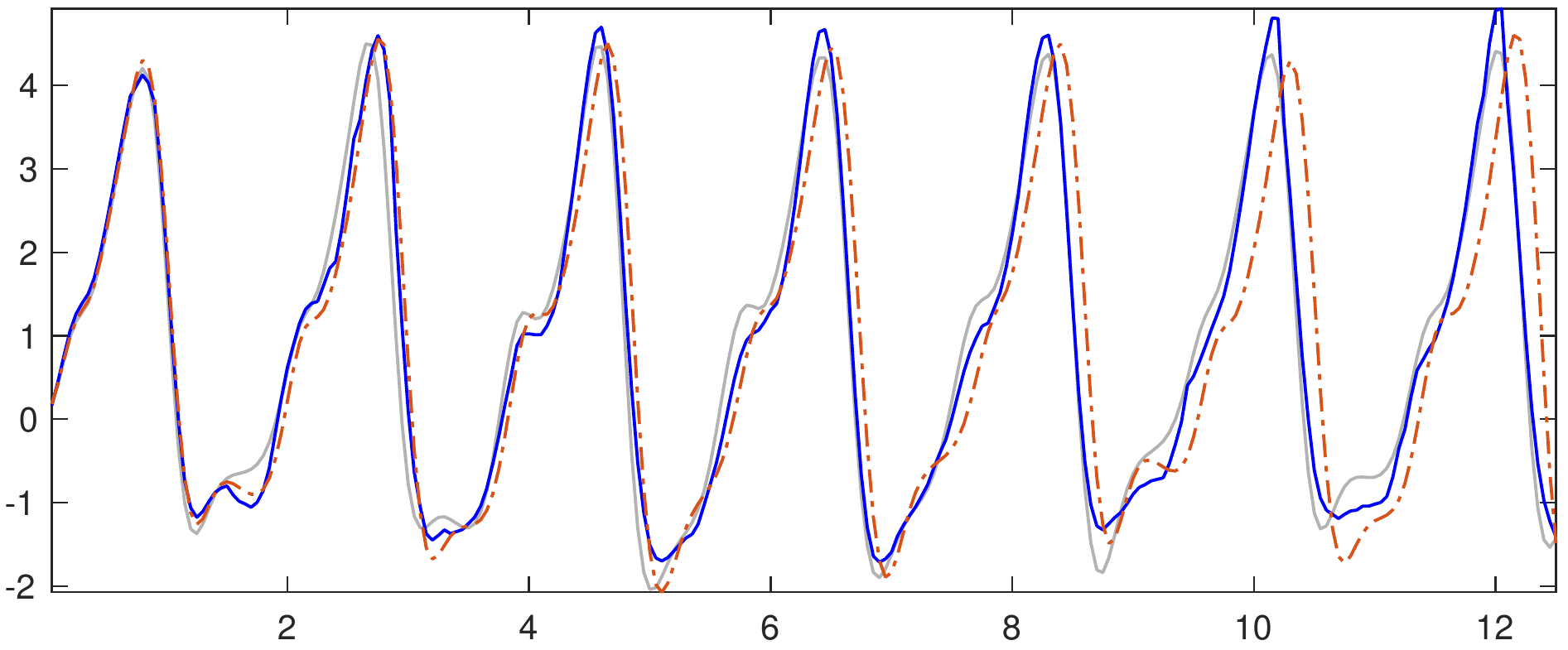}
\put(8,37) {\footnotesize Quasiperiodic} 
\end{overpic}

\begin{overpic}[width=.6\textwidth]{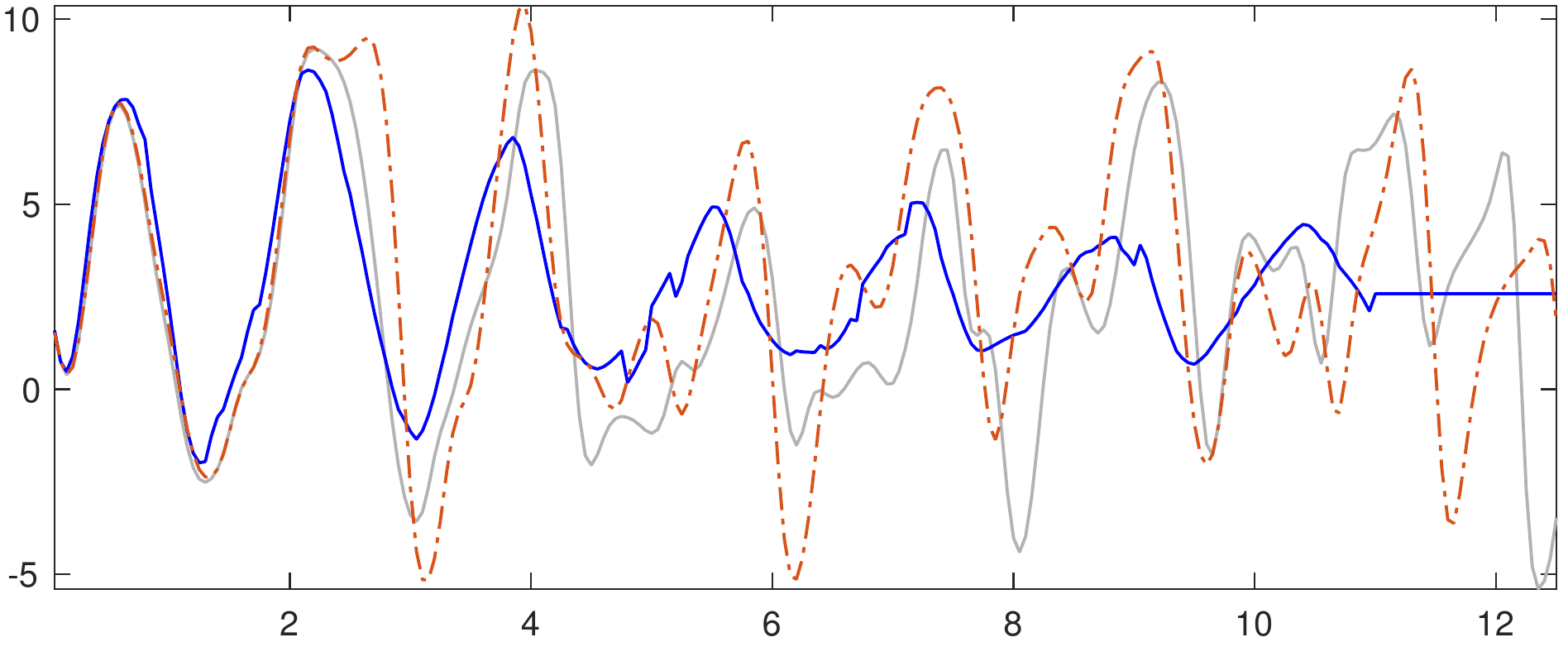} 
\put(50,-2) {$\tau$}
\put(8,37) {\footnotesize Chaotic} 
\end{overpic}
\caption{{\bf Comparison of data-driven and model-data-driven prediction} 
The true trajectory is shown in grey, the KAF data-driven prediction
in blue and the model-data-driven predictions based on \eqref{eq:l96b} in dotted-red;
the periodic, quasiperiodic, and chaotic regimes are considered in turn.}
\label{fig:comparegp}
\end{figure}

Explicit details of the procedure we use to build a GP closure are 
described in Appendix~\ref{as:GPclosure};
here we observe that for training we use tuples
$\left\{ x_k(t_n), (By)_k(t_n)  \right\}_{n=1}^N$, over all $k=1,\cdots, 9.$
See Figure~\ref{fig:gpr} to see the data used (red random subsamples, without
replacement, of the total grey data set), and an approximate
GP closure $\cgp$ determined from that data.

\begin{figure}[!ht]
\centering
\begin{overpic}[width=.6\textwidth]{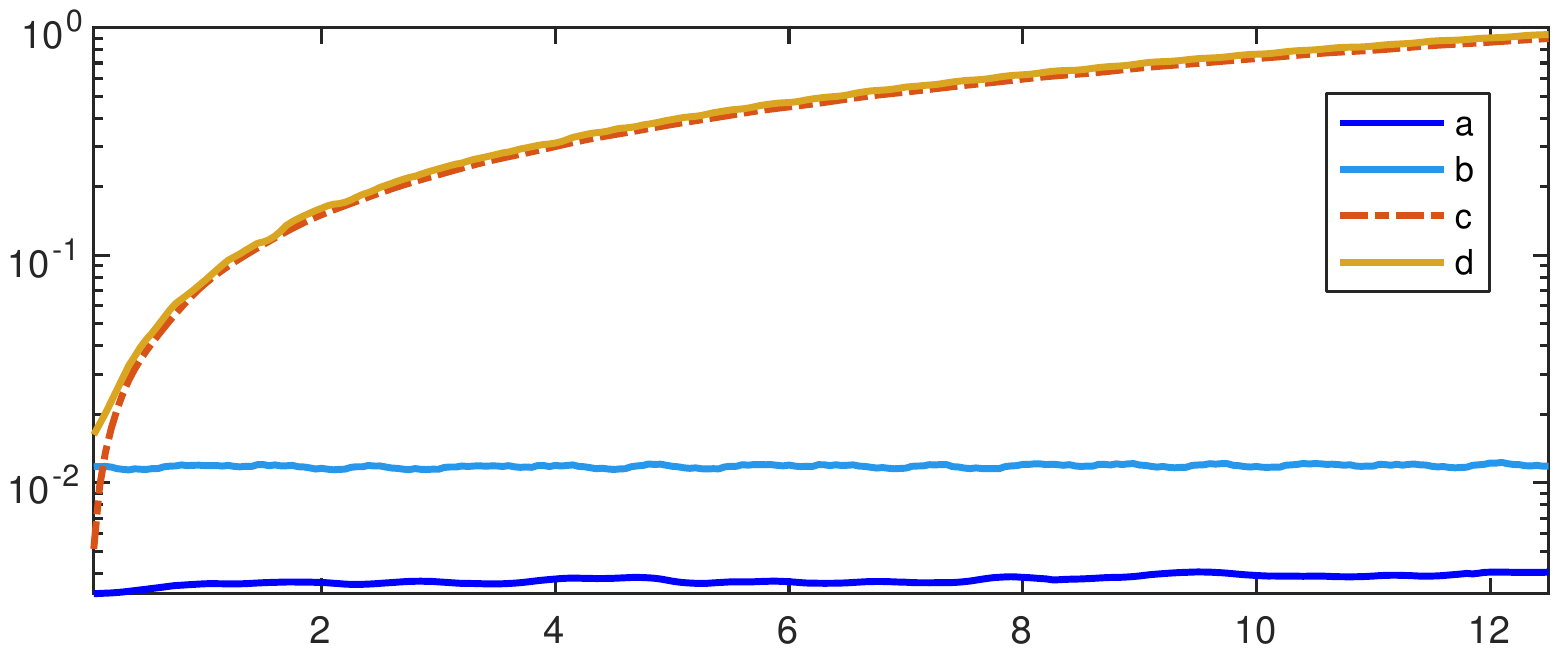} 
\put(8,37) {\footnotesize Periodic} 
\end{overpic}%

\begin{overpic}[width=.6\textwidth]{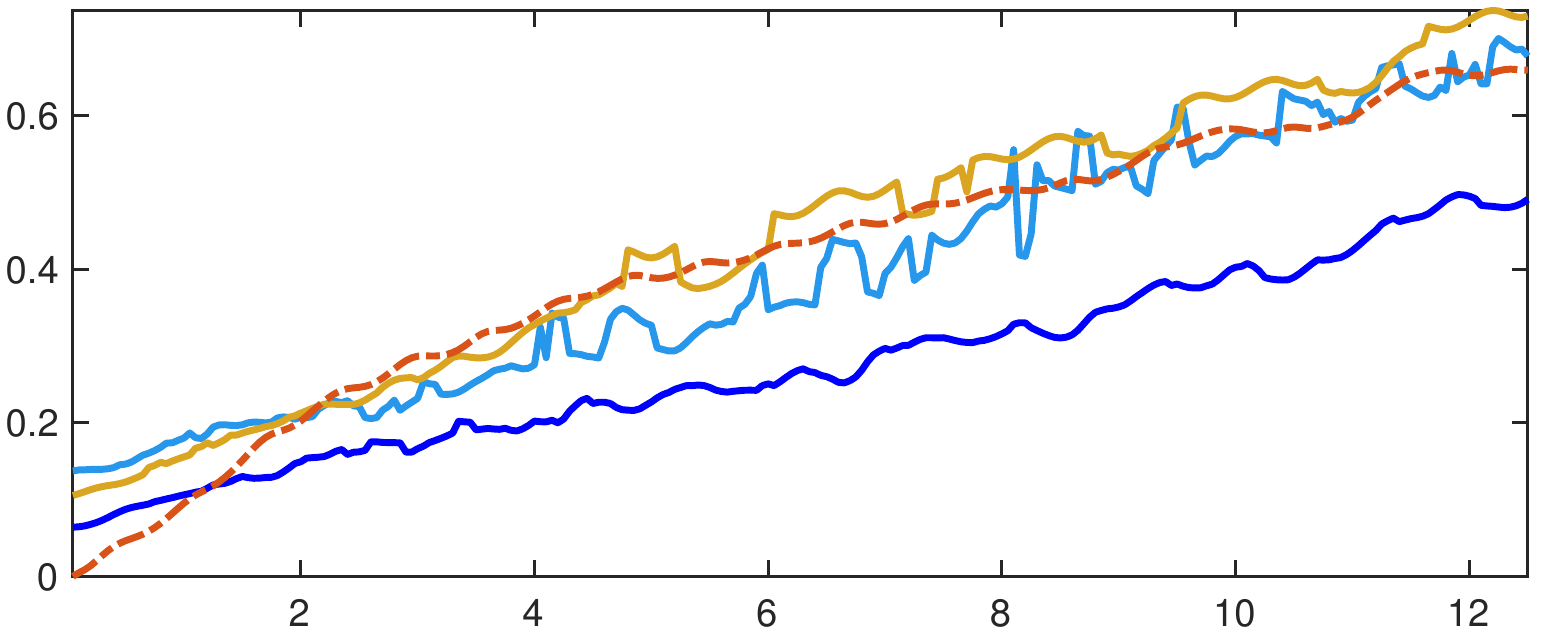}
\put(8,37) {\footnotesize Quasiperiodic} 
\end{overpic}%

\begin{overpic}[width=.6\textwidth]{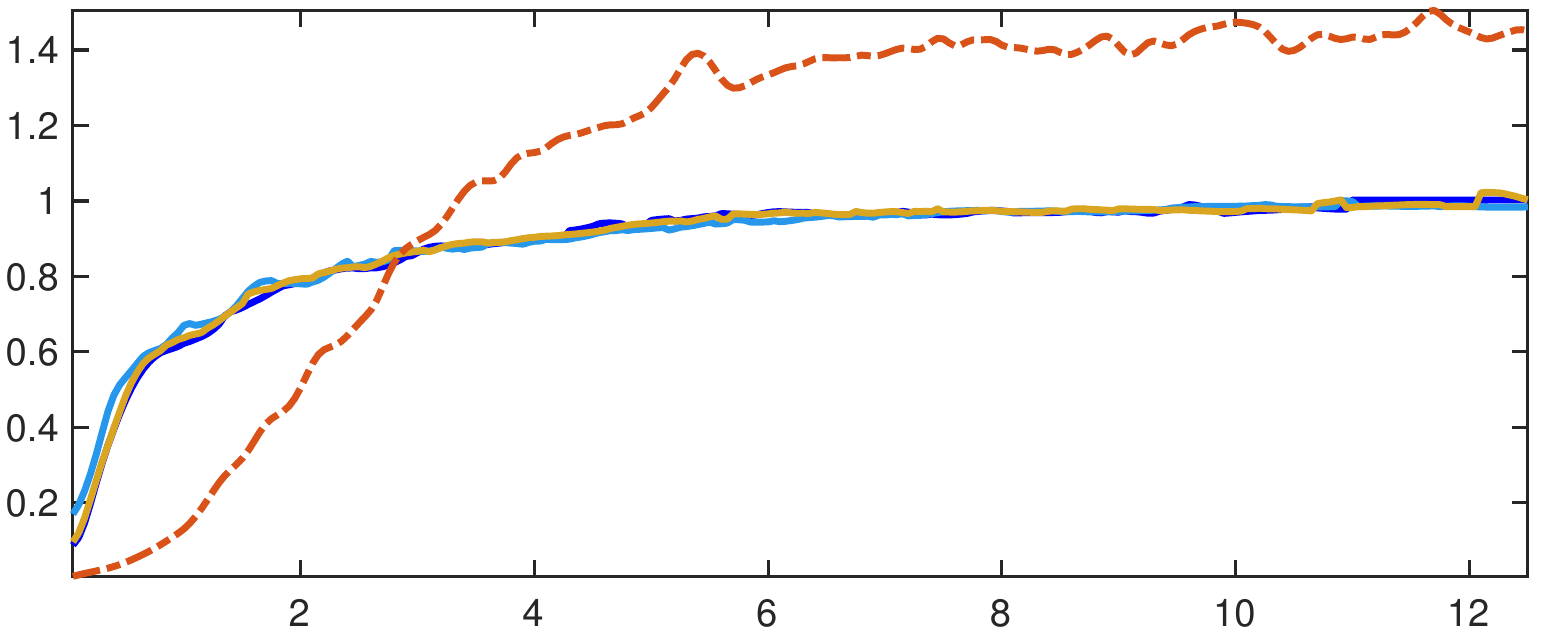} 
\put(50,-2) {$\tau$}
\put(8,37) {\footnotesize Chaotic} 
\end{overpic}%
\caption{{\bf \RMSE comparison for the four cases a)--d) described in the
text, in the periodic, quasiperiodic, and chaotic regimes.} In the periodic regime, KAF (a) is an ideal predictor with negligible growth in error (note the logarithmic scale). In the quasiperiodic response regime, the growth in \RMSE with KAF (a,b) is significantly slower than that of the GP-based
ODE prediction (c). In the chaotic response regime, the GP-based
ODE prediction (c) is more accurate in the near term, 
yet KAF error stabilizes as the prediction converges to the conditional mean.}
\label{fig:rmse_abcd}
\end{figure}

Once we have the closed model appearing in \eqref{eq:l96b} we may use
it to predict the variable $x$ appearing in \eqref{eq:l96}, and we
may compare that prediction with the one made by KAF.  
Figure~\ref{fig:comparegp} shows the result of doing so. It shows
that the KAF approach is superior in the periodic and quasi-periodic
settings, but that for predictions of the trajectory itself the
model-data based predictor \eqref{eq:l96b} is superior to KAF in the 
chaotic case. Note that the model-data based predictor has access
to more data than does the KAF, and requires model knowledge; the
KAF is entirely data-driven.

We now dig a little deeper into the comparison.
We do this in a systematic way in the
periodic, quasi-periodic and chaotic regimes. In each of these
three cases we show four \RMSE error curves, labelled as follows: 
a) the standard KAF based on $x$ data alone; b) an enhanced KAF 
using $(x,By)$ data, the same data used to train the ODE \eqref{eq:l96b}; 
c) a prediction using the ODE \eqref{eq:l96b}; 
d) a KAF prediction trained on $X$ data alone, generated by the 
ODE \eqref{eq:l96b}. Figure \ref{fig:rmse_abcd} shows
that KAF a) is the ideal predictor in the periodic regime
and is  near-ideal in the quasi-periodic regime; on the other
hand the ODE \eqref{eq:l96b} predictor c) is ideal for short-term 
predictability in the chaotic case. 
Augmenting observations with $By$ within KAF, as in b), gives errors
similar to those arising from a), when observing $x$ alone; thus knowledge
of $By$ provides little extra information. 
In the chaotic case, the {\RMSE}s of KAF trained on $x$, a), and 
on $X$, d), are very close, confirming that the ODE \eqref{eq:l96b}
for $X$ captures the invariant measure of the approximately Markovian
variables $x$ as intended.

\revisionOne{We emphasize the difference between averaging and homogenization
here: in the averaging case observing the fast variables adds nothing to our
prediction because there is a closed system determined only by the slow
variables (see Figure~\ref{fig:rmse_abcd}, graphs a) and b) in all three plots).
By contrast, in the homogenization case, observing the fast variables improves
short-term predictions because it provides further information about the driving
stochastic process entering the homogenized limit (see
Figure~\ref{fig:nonMarkovL63}).}

\begin{figure}[!ht]
\centering
\begin{overpic}[width=.6\textwidth]{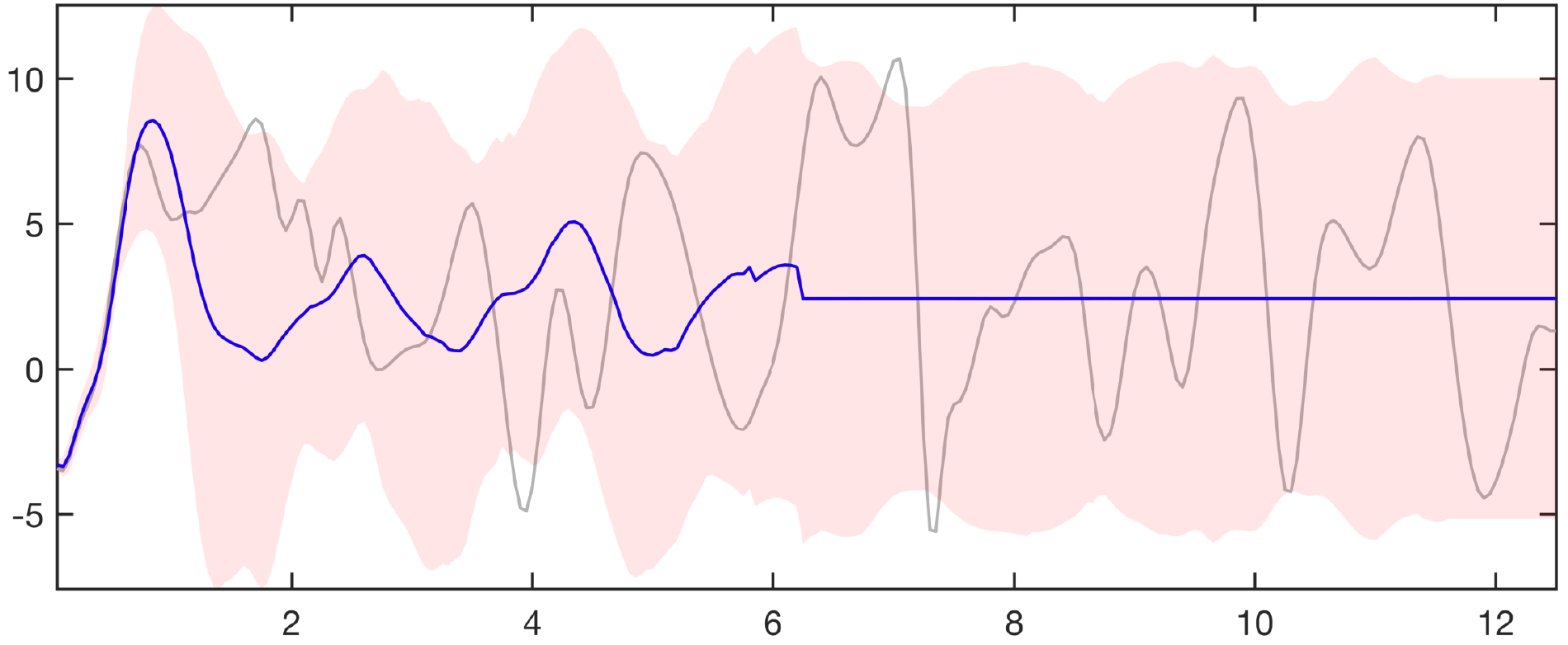} 
\put(40,37) {\footnotesize $\ix(\omega) = x$} 
\put(100,33){\fbox{%
  \parbox{2.3cm}{\scriptsize
    {\color{gray} \bf --}  Trajectory \\
    {\color{blue} \bf --}  Prediction\\
    {\color{pink}$\blacksquare$} Predicted 2$\sigma$ 
  }%
}}
\end{overpic}
\caption{{\bf Prediction in non-Markovian regime, $F_x=10,~\lorenzeps= 1$, }results in much shorter accurate trajectory predictability, 
followed by rapid convergence of the conditional mean to a constant. Note,
however, that the uncertainty prediction bands contain the true trajectory 
for all time.}
\label{fig:nonMarkovL96}
\end{figure}

\begin{figure}[!ht]
\centering
\begin{overpic}[width=.35\textwidth]{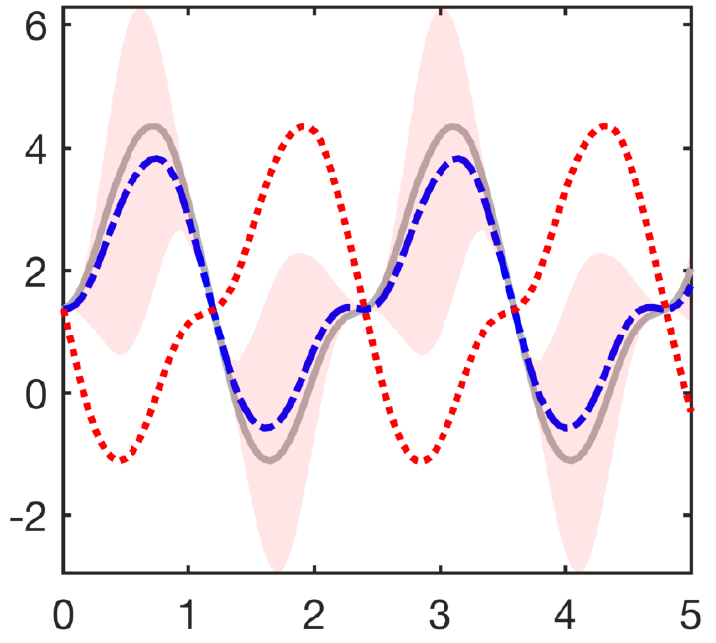} 
\put(30,82) {\footnotesize $x_1(0) = 1.3736$} 
\put(48,-2) {$\tau$} 
\end{overpic}
\begin{overpic}[width=.35\textwidth]{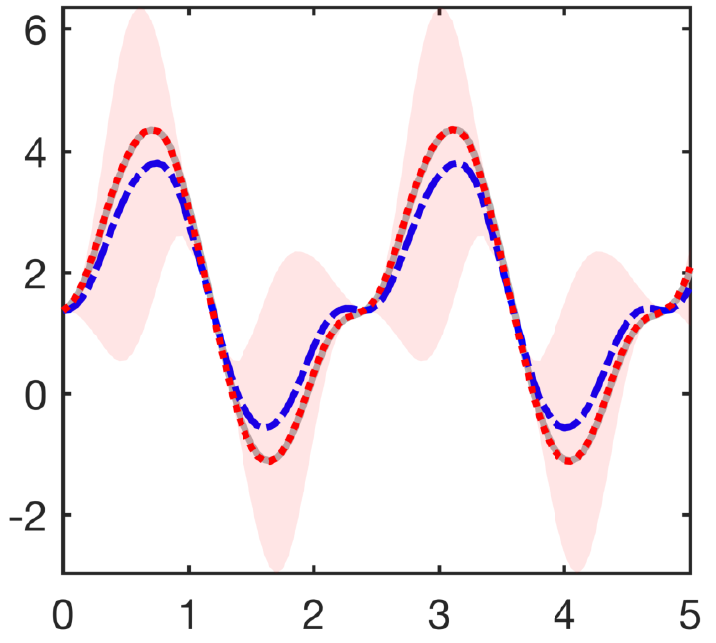} 
\put(30,82) {\footnotesize $x_1(0) = 1.3799$} 
\put(-15,92) {\fbox{\footnotesize $\ix(\omega) = x_1$}}
\put(48,-2) {$\tau$} 
\put(90,73){\fbox{%
  \parbox{2.55cm}{\scriptsize
    {\color{gray} \bf --}  Trajectory \\
    {\color{red} \bf ...}  Lorenz' method \\
    {\color{blue} \bf - -}  KAF prediction \\
    {\color{pink}$\blacksquare$} KAF Predicted 2$\sigma$
  }%
}}
\end{overpic}
\caption{\revisionOne{{\bf Lorenz' method~\eqref{eq:LAF} sensitivity to initial conditions,} in the partially-observed periodic $F_x=5, \varepsilon^{-1}=128$ regime, results in diverging predicted trajectories for nearby initial conditions, separated only by 0.006. By contrast, KAF, which is continuous with respect to initial condition, shows moderate predictive skill and makes nearly identical predictions for nearby initial conditions.}}
\label{fig:LAFcomparison}
\end{figure}
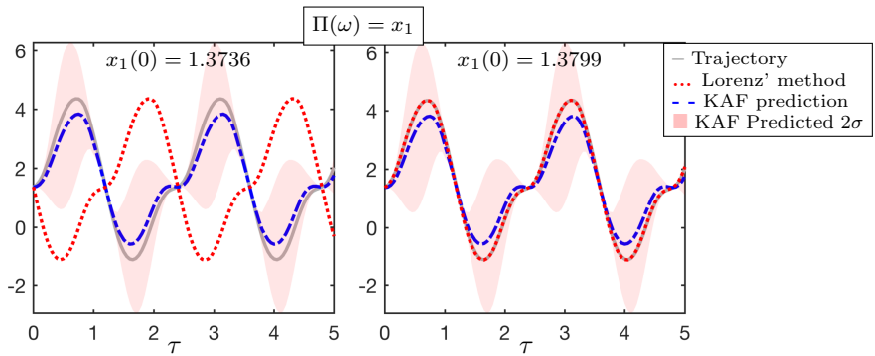

\subsection{Non-Markovian regime}
\label{ssec:NM2}

In the preceding subsections we studied predictors for $x$, based only
on time-series data in the $x$ coordinate, for the equation \eqref{eq:l96}.
We studied the scale-separated regime where $\lorenzeps \ll 1$ and $x$
is approximately Markovian and deterministic -- it is approximately
governed by an ODE. Here we study the behavior of identical
predictors when $\lorenzeps=1$; the system ~\eqref{eq:l96} then no longer
exhibits averaging and $x$ is no longer Markovian because there
is no scale-separation between $x$ and $y$.
This experiment is conducted with $F_x=10$. 
Because of the lack of Markovian behaviour we expect rapid
loss of predictability in time, when $\ix(\omega)=x, F(\omega) = x_1$. 
The resulting conditional mean and variance, shown in 
Figure~\ref{fig:nonMarkovL96}, confirms this intuition. 
Indeed the conditional mean is out of phase with the truth
at lead time $\tau=1$, and this is also reflected in the 
large growth of the conditional variance. Furthermore, the 
conditional mean tapers to a constant at $\tau=6$, twice as 
quickly as it does in the $\lorenzeps \ll 1$ setting in which 
this tapering occurs at $\tau\approx 11$ 
(Figures~\ref{fig:forecast96},\ref{fig:comparegp}).

%

\revisionOne{\subsection{Comparison with Lorenz' method} We illustrate the advantage of KAF over Lorenz' original method of analog forecasting~\eqref{eq:LAF}, which can produce predictions that are discontinuous with respect to initial condition. In particular, this occurs when data are partially observed from a larger state space, and different states map to identical partial observations. 
To study this, we observe a single coordinate $x_1$ of the periodic regime ($F_x = 5, \varepsilon^{-1} = 128 $) so that  the observed data are highly non-Markovian. We select initial conditions that are $\mathcal{O}(10^{-3})$ apart, but are separated in time by integer multiples of the period.   Figure~\ref{fig:LAFcomparison} plots the resulting predictions from Lorenz' method and KAF. Although Lorenz' method is accurate for one initial condition (right), it gives a diverging prediction for a nearly identical point (left). By  contrast, KAF is continuous with respect to initial condition and displays theoretically optimal predictive skill (in an RMSE sense) for even highly non-Markovian observation data. This experiment illuminates a key feature of KAF: that it gives consistent predictions that are continuous with respect to initial conditions. Note, also, that KAF uncertainty predictions of a periodic observable are also periodic, and vanish at every half period when predictions intersect the ground truth.}

\section{Conclusions}
\label{sec:C}

\begin{enumerate}

\item We have studied KAF for data-driven prediction:

\begin{itemize}
\item we use multi-scale systems to create dynamical systems
in which a subset of the variables (\emph{slow variables}) 
evolve in an approximately Markovian fashion;

\item we study KAF performance for a range of systems in which the slow
variables are governed approximately by stochastic, chaotic,
quasi-periodic and periodic behaviour;
 
\item in the stochastic case we use the homogenized equations
for the slow variables to obtain explicit formulae for the
eigenfunctions of the operator underlying KAF, and use these
to validate the performance of the KAF method; 

\item in the chaotic, quasi-periodic and periodic cases we use 
the averaged equations for the slow variables to obtain a GPR-based
approximate closure model, and use hybrid data-model predictions, 
from this closure, in order to evaluate the KAF method.

\end{itemize}

\item What we illustrate about use of the KAF:

\begin{itemize}

\item when the variable being predicted is (approximately) a
component of a stochastic Markov process then prediction of individual 
trajectories is not possible, whilst the mean and variance, averaged over
possible realizations of the stochastic behaviour, can be accurately predicted
by KAF; 

\item when the variable being predicted is (approximately) a
component of a deterministic but chaotic Markov process then prediction 
of individual trajectories is also not possible, except over short time
horizons; 

\item 
    in both the chaotic and stochastic cases, whilst prediction of 
individual trajectories is not to be expected, simply bounding
the future trajectory may be useful in applications -- to this 
end we show that in all cases two standard deviation bands around
the predicted conditional mean always reliably capture the truth;

\item in both the stochastic and chaotic settings a signature of the
lack of predictability is the convergence of the predictor to a constant,
for large $\tau$, accompanied by the data-driven choice of parameter $\ell$
converging to one;

\item when the variable being predicted is (approximately) a
component of a deterministic quasi-periodic or periodic process
the prediction of individual trajectories over long time horizons
is possible; in this case parameter $\ell$ stays away from one,
for significantly large $\tau;$ 

\item in all cases the predicted standard deviations around the
predictor provide a reliable indicator of the time-scale on which
the predictor is accurate trajectory-wise, and on longer time-scales
provide a reliable indicator of the scale of the errors incurred,
trajectory-wise;

\item the choice of kernel, and the interpretation of eigenfunctions as 
harmonics over the observed submanifold, is well-adapted to capturing 
periodic or quasi-periodic structure whilst the design of eigenfunctions
to include the constant eigenfunction permits convergence to the 
conditional mean for chaotic, mixing dynamics;

\item the choice of $\ell$ indicates the number of harmonics 
contained in the predicted variable and is often approximately
constant in $\tau$ in periodic and quasi-periodic settings,
whilst it converges to one in the chaotic or stochastic mixing cases;

\item when the variable being observed does not evolve in an approximately
Markovian fashion then the KAF conditional mean cannot track $x$ for even 
short times, as observed for $\lorenzeps = 1$ for both the
Lorenz 63 and Lorenz 96 examples.

\end{itemize}

\item The work also suggests a number of directions for future study
in the area of KAF: 

\begin{itemize}
\item delay embedding can be used to deal with non-Markovian
behaviour and it would be of interest to automate the choice of
delay embedding dimension within the KAF framework, to get
closer to Markovianity;
\item combining KAF with data assimilation holds the possibility
of greater predictability -- for work in this direction see~\cite{hamilton2016ensemble} which studies the EnKF with a data-driven 
model update; 
\item it would be of interest to have algorithms to identify slow 
subspaces, using data living in larger spaces, and hence to conduct
KAF using approximately Markovian variables;
\item it would be of interest to extend KAF predictor $Z_{\tau}$ so that
it acts on the joint space of initial conditions and key parameters, enabling
prediction at as yet unseen parameters.
\end{itemize}
\end{enumerate}

\vspace{0.1in}

\noindent{\bf Acknowledgments:} The work of DB and AMS is supported
by the generosity of Eric and Wendy Schmidt by recommendation of 
the Schmidt Futures program, by Earthrise Alliance, Mountain Philanthropies, 
the Paul G. Allen Family Foundation, and the National Science Foundation 
(NSF, award AGS1835860). AMS is also supported
by NSF (award DMS-1818977) and by the Office of Naval Research 
(award N00014-17-1-2079). DG is grateful to the Department of Computing and Mathematical Sciences at the California Institute of Technology for hospitality and for providing a stimulating environment during a sabbatical, where part of this work was completed. DG is supported by NSF (awards 1842538 and DMS-1854383) and ONR (awards N00014-16-1-2649 and N00014-19-1-242). KM is supported by the NSF Mathematical Sciences Postdoctoral Research Fellowship (award 1803663).

\vspace{0.1in}

\newpage
\bibliographystyle{siamplain}
\bibliography{references,references_dg}

\appendix
\section{\label{appComputation}Computation}
\label{a:a}

This appendix contains details of the implementation of KAF. Subsection~\ref{as:la} describes the algorithms for computing the kernel, diffusion eigenbasis, RKHS basis functions and finally, the prediction. Our specific choice of 
kernel, which  endows the RKHS structure, is explained and motivated
in subsection~\ref{as:kernel}. We outline procedures for choosing truncation parameter $\ell$ and approximating the conditional variance of the forecast in subsections~\ref{as:L} and~\ref{as:V}.

\subsection{Linear Algebra}

\label{as:la}


\subsubsection{Algorithm 1 (Diffusion eigenbasis)}
\label{sssec:Alg1}
The starting point of this algorithm is the variable-bandwidth diffusion kernel function $\est{\kappa}:\mathcal{X}\times \mathcal{X}\to \R$ 
\begin{displaymath}
\est{\kappa}(x,x')=\exp\left(\frac{-|x-x'|^2}{\kerneleps\est{r}(x)\est{r}(x')}\right),
\end{displaymath}
described further in Appendix~\ref{as:kernel}. 
An automatic procedure for estimating the data-dependent bandwidth function $\est{r}$ and width $\kerneleps$ are given in~\cite{berry2015nonparametric}.

\begin{itemize}
\item Inputs
\begin{itemize}
\item Training data $x_0,x_1,\dots x_{N-1}\in\mathcal{X}$ 
\item Desired maximum number of eigenvectors $L$
\end{itemize}
\item Outputs
\begin{itemize}
\item Diffusion eigenvectors $\phi_0,\dots, \phi_{L-1}\in\R^N$, stacked in matrix $\bPhi\in\R^{N\times L}$
\item Diffusion eigenvalues $\lambda_0,\dots,\lambda_{L-1}\in\R$ in diagonal matrix $\bLambda\in\R^{L\times L}$
\end{itemize}
\item Steps
\begin{enumerate}
\item Form the matrix $\bK\in\R^{N\times N}$ with entries $K_{ij} =\kappa_N(x_i,x_j)/N$.
\item Compute $\est{v},\est{w}\in\R^N$ using $\est{v}=\bK\ones$ and $\est{w} = \bK V^{-1}\ones$, where $V = \mbox{diag}(\est{v})$.
\item Form normalized kernel matrix $S = V^{-1}\bK W^{-1/2}$ where  $W = \mbox{diag}(\est{w})$.
\item Compute $L$ largest singular values $\sigma_0,\dots,\sigma_{L-1}$ and corresponding left singular vectors $\phi_0,\dots,\phi_{L-1}$ of $S$. Stack eigenvectors columnwise into $N\times L$ matrix $\bPhi := [\phi_0,\dots,\phi_{L-1}]$ and form $L\times L$ diagonal matrix of eigenvalues $\bLambda := \mbox{diag}(\lambda_j)$, with $\lambda_j := \sigma_j^2$. 
\end{enumerate}
\end{itemize}

\paragraph{Remark} Recall that the key idea of KAF is the eigendecomposition of the Markov operator $G$ which can be represented by an $N\times N$ matrix
$$  G\phi = \lambda \phi.$$
The matrix $G$ is never explicitly formed, instead we exploit the fact that $G=SS^T$, and hence $\phi$ are the left singular vectors of $S$. Thus we bypass the eigendecomposition step with a reduced singular value decomposition (SVD) in step (4). Note that the SVD approach is natural when working with kernels with an explicit factorization $ G = SS^T$, including the bistochastic kernels described in Appendix~\ref{as:kernel}. For more general kernels one typically performs direct eigendecomposition of $G$. KAF can also be implemented with non-symmetric kernels satisfying a detailed balance condition making them conjugate to positive-definite kernels; see Section~4.2 in \cite{AlexanderGiannakis20} for more details.

\subsubsection{Algorithm 2 (RKHS basis functions) }
\label{sssec:Alg2}

RKHS basis functions are computed using the Nystr\"om extension~\eqref{eq:psiN} reproduced below
\begin{displaymath}
  \psi_j = \frac{1}{N \lambda_j^{1/2}} \sum_{n=0}^{N-1} k( \cdot, x_n )\phi_j(x_n), \quad \lambda_j > 0.
\end{displaymath}

\begin{itemize}
\item Inputs
\begin{itemize}
\item Out-of-sample data $\xo_0,\xo_1,\dots \xo_{\hat{N}-1}\in\mathcal{X}$ 
\item Diffusion eigenvectors $\bPhi\in\R^{N\times L}$ 
\item Diffusion eigenvalues $\bLambda \in \R^{L\times L} $
\end{itemize}
\item Outputs
\begin{itemize}
\item RKHS basis $\bPsi = [\psi_0,\dots, \psi_{L-1}]$ 
\end{itemize}
\item Steps
\begin{enumerate}
\item Form matrix $\hat{S}\in\R^{L\times N}$ with entries $\hat{S}_{ij} = \est{\kappa}(\xo_i,x_j)$. Note that this requires another kernel density estimation step to evaluate the sampling density on out-of-sample data.
\item Compute the RKHS basis matrix $\bPsi = \bLambda^{-1/2}\hat{S}\bPhi/N$. RKHS basis functions (evaluated at the out-of-sample points) are the columns $\psi_0,\dots, \psi_{L-1}$ of $\bPsi$.
\end{enumerate}
\end{itemize}

\subsubsection{ Algorithm 3 (Predictor)}
\label{sssec:Alg3}

The final predictor is constructed according to~\eqref{eq:Z_EXPAND}, reproduced here for convenience
\begin{displaymath}
Z_\tau(x) = \frac{1}{N}\sum_{n=0}^{N-1} \left(\sum_{j=0}^{\ell(\tau)-1} \frac{ \psi_j(x)\phi_j(x_n)}{\lambda_j^{1/2}} \right) f_{n+q}.
\end{displaymath}

\begin{itemize}
\item Inputs
\begin{itemize}
\item Lead time $\tau = q\Delta t$
\item Truncation parameter $\ell \le L$
\item Vector of sampled observable $\ff_\tau = [f_{q}, \dots, f_{N-1+q}]^T$ 
\item Diffusion eigenvectors $\bPhi_{\ell}\in\R^{N\times\ell}$ (first $\ell$ columns of $\bPhi$)
\item Diffusion eigenvalues $\bLambda_{\ell}=\mbox{diag}(\lambda_0,\dots,\lambda_{\ell-1})$ 
\item RKHS basis $\bPsi_{\ell}\in\R^{N\times\ell}$ (first $\ell$ columns of $\bPsi$)
\end{itemize}
\item Output
\begin{itemize}
\item Prediction $Z_\tau$ at $\xo_0,\xo_1,\dots \xo_{\hat{N}-1}$ 
\end{itemize}
\item Steps
\begin{enumerate}
\item Form $\ell$-dimensional coefficient vector $\bc := \bPhi_{\ell}^T \ff_\tau /N$
\item Compute $\bz := \bPsi_{\ell}\bLambda_{\ell}^{-1/2}\bc$. 
Report prediction for $i^{th}$ initial point $Z_\tau(\xo_i) := z_{i}$. 
\end{enumerate}
\end{itemize}

\subsection{Choice of Kernel}
\label{as:kernel}
Details concerning the kernel choice may be found in section 5 of
\cite{AlexanderGiannakis20}. Here we briefly summarize the key ideas.
Our starting point is the Gaussian kernel
$$\kappa(x,x';\kerneleps) = \exp\bigl(-|x-x'|^2/\kerneleps\bigr)$$
where we assume that $\X$ is a subset of $\R^d$, $x, x' \in \X$
and $|\cdot|$ denotes the Euclidean norm on $\R^d.$ 
Note that this kernel is data-independent.
From it we can generalize to a data-dependent kernel defined by \cite{BerryHarlim16}
\begin{align*}
\est{\kappa}(x,x')&=\exp\bigl(-|x-x'|^2/(\est{r}(x)\est{r}(x')\kerneleps)\bigr),\\
\est{r}(x)&= \est{q}(x)^{\frac{1}{m}},\\
\est{q}(x)&=\frac{1}{(\pi \delta)^{m/2}}\int_{\Omega} \kappa(x,\ix(\omega);\delta)\est{\mu}(d\omega). 
\end{align*}  
Here, $\kerneleps, \delta$ and $m$ are lengthscale and dimension parameters, respectively, estimated from the data. The role of including
$r_n$, and hence a variable bandwidth kernel, is to compensate
for variations in sampling density across the space. 
A key conceptual idea underlying the construction of $\est{q}$ is that
it approximates the Lebesgue density of the measure $\mu$ in the large
data limit. Thus the kernel $\est{\kappa}$ weights the distance of
$x$ and $x'$ in a manner which reflects the sampling density of the
data.

From $\est{\kappa}$ the kernel $\est{k}$ is constructed as follows, 
invoking a second principle which is to design a 
bistochastic Markov kernel \cite{CoifmanHirn13}.
Doing so ensures that the top eigenvalue of $G$ is $1$ with
corresponding eigenvector a constant; then the hypothesis space
also contains constants. This is useful for capturing the mean
of the predicted quantity $U^{\tau}F$ and, in particular, plays
a central role in the large $\tau$ asymptotics for mixing systems.
To achieve the Markov property we proceed as follows. First we define
\begin{align*}
\est{v}(x)&=\int_{\Omega} \est{\kappa}\bigl(x,\ix(\omega)\bigr)\est{\mu}(d\omega),\\
\est{w}(x)&=\int_{\Omega} \frac{\est{\kappa}\bigl(x,\ix(\omega)\bigr)}{\est{v}\bigl(\ix(\omega)\bigr)}\est{\mu}(d\omega ).  \end{align*}
Since the above empirically determined quantities only take values at the $N$ sampled points, they are isomorphic to vectors in $\R^N$ identified by the same name within Algorithm~\ref{sssec:Alg1}. Finally define
$$\est{k}(x,x')=\int_{\Omega}\frac{\est{\kappa}\bigl(x,\ix(\omega)\bigr)\est{\kappa}\bigl(\ix(\omega),x'\bigr)}{\est{v}(x)\est{w}\bigl(\ix(\omega))\est{v}(x')} \est{\mu}(d\omega).$$
The operators constructed from the unnormalized and normalized kernels, $\est{\kappa}$ and $\est{k}$, are likewise represented by $N\times N$ matrices $K$ and $S$, respectively, in Algorithm~\ref{sssec:Alg1}. 
It may be verified that the Markov property is satisfied and so too are
the positivity conditions required for the aforementioned large data
convergence result.

\subsection{Choice of $\ell$}
\label{as:L}

Recall that the predictor \eqref{eq:Z_PRED} actually corresponds to a 
{\em family} of predictors parameterized by the desired lead time $\tau$ and
by the truncation parameter $\ell$. Thus we write $Z_{\tau,\ell}.$ Here we 
describe how to choose $\ell$ for a fixed lead time $\tau$. In practice, 
the choice of $\ell$ is determined from the minimizer of the empirical 
\RMSE based on~\eqref{eq:LN}, computed from a validation data set 
with $\tilde{N}$ samples $\xo_0,\dots,\xo_{\tilde{N}-1}:$ 
\begin{displaymath}
\ell = \argmin_{\ell' =1,\dots,L} \RMSE(Z_{\tau,\ell'}).
\end{displaymath}

\subsubsection{Algorithm 4 (Tuning $\ell$)}
\label{sssec:tuning}

\begin{itemize}
\item Inputs
\begin{itemize}
\item Forecast lead time $\tau = q\Delta t$
\item Validation out-of-sample data $\xo_0,\dots,\xo_{\tilde{N}-1}\in \mathcal{X}$
\item Ground truth vector of observables $\hat{\ff}_\tau =[\fo_q,\dots,\fo_{\tilde{N}-1+q}]^T$
\item Diffusion eigenvectors $\bPhi$ from Algorithm 1
\item Diffusion eigenvectors $\bLambda$ from Algorithm 1
\end{itemize}
\item Outputs
\begin{itemize}
\item Truncation parameter $\ell$
\end{itemize}
\item Steps
\begin{enumerate}
\item Compute RKHS basis functions $\bPsi_L$ using Algorithm 2. Set $R := \infty$.
\item for $\ell' = 1$ to $L$
\begin{enumerate}
\item Compute predictor $Z_{\tau,\ell'}(\xo):=\bz$ using Algorithm 3.
\item Compute $\RMSE_{l'}$ for $Z_{\tau,\ell'}(\xo)$ as $\RMSE_{\ell'} := \| \bz - \hat{\ff}_\tau \|_2$.
\item if $\RMSE_{\ell'}\le R$ set $\ell=\ell'$
\end{enumerate}
\item Return $\ell$.
\end{enumerate}
\end{itemize}
This tuning procedure must be carried out for each desired lead time $\tau$. 

\subsubsection{Long-time behavior of $\ell$}
\label{sssec:longt}

It should be noted that in the presence of mixing or chaotic dynamics, for long lead times $\tau$ the projected subspaces become one-dimensional, i.e., $\ell=1$ and the predictor converges to a constant (this occurs only if the subspace includes the constant eigenfunction from a Markov kernel operator). The weak convergence of the conditional expectation of mixing dynamics to the mean of the observable $F$ is described in~\cite{AlexanderGiannakis20},
\begin{equation}
\lim_{\tau\rightarrow\infty} \langle g, \Expect[U^\tau F | \ix] - \Expect[F]\mathbf{1}_\Omega \rangle_{L^2_\mu} = 0, 
\end{equation}
a property that is a direct consequence of a measure-theoretic definition of mixing
\begin{equation}
\lim_{\tau\rightarrow\infty} \langle U^{\tau *}g, h\rangle_{L^2_\mu} = \Expect[g]\Expect[h], 
\end{equation}
where the expectations are taken over the invariant measure $\mu$.

\subsection{Formula for Variance}
\label{as:V}

The uncertainty associated with the prediction at each lead time can be estimated via the conditional variance, 
\begin{equation*}
\var[U^\tau F| \ix] = \Expect[(U^\tau F - \Expect[U^\tau F|
\ix])^2 | \ix] \approx \Expect[(U^\tau F - Z_\tau)^2 | \ix].
\end{equation*}
The variance is yet another observable in $\ltm$ that can be evaluated using the same basis functions as the predictor, and, once the predictor is computed, it only remains to compute the expansion coefficients
\begin{equation}
             \hat{c}_j(\tau) = \langle \phi_j \circ \ix, (U^\tau F- Z_\tau )^2 \rangle_{L^2(\est{\mu})} = \frac{1}{N} \sum_{n=0}^{N-1} \phi_j(x_n)(f_{n+q}-Z_\tau(x_n))^2.
\end{equation}
Note that this computation requires a different optimal truncation $\ell(\tau)$ for the variance, and hence, another validation set for parameter tuning:
\begin{enumerate}
\item (Tuning $\ell$) Run Algorithm 4 on a separate validation data set for which the predictor is already computed, using the new observable $\hat{g}_\tau$ instead of $\hat{f}_\tau$
\begin{displaymath}
\hat{g}_\tau = [(\fo_q-Z_\tau(\xo_0))^2,\dots,(\fo_{\tilde{N}-1+q}-Z_\tau(\xo_{\tilde{N}-1}))^2]^T.
\end{displaymath}
\item Run Algorithm 3 on the initial prediction data using the observable $g_\tau$ instead of $f_\tau$
\begin{displaymath}
g_\tau = [(f_q-Z_\tau(x_0))^2,\dots,(f_{N-1+q}-Z_\tau(x_{N-1}))^2]^T,
\end{displaymath}
denoting the output predicted variance as $V_\tau$. 
\item Finally, the uncertainty bands of the predictor at each lead time are given by two standard deviations, or twice the square root of the variance, $Z_\tau\pm 2\sqrt{|V_\tau|}$.
\end{enumerate}



\section{Details Of The GP Closure}
\label{as:GPclosure}

In this section we describe details of the construction of
the Gaussian process (GP) underlying the model-data-driven 
approach, and leading to the approximate closed equation
\eqref{eq:l96b}. We use \LNS{} explicitly, but the methodology 
easily generalizes to other multiscale systems.

Construction of a GP closure is performed using the following steps:
\begin{enumerate}[(a)]
  \item choose random initial conditions for \LNS{};
  \item numerically integrate for time $T_{\text{conv}}$ to determine
an initial condition on the global attractor;
  \item numerically integrate from this initial condition
with a fixed time-step $\Delta t$ for time
    $T_{\text{learn}}$ to collect pairs of data:
    \[
      \left\{ \bigg( x_k(t_n),
      (By)_k(t_n) \bigg) \right\}_{n=1}^N,
    \]
for $k=1,\dots, 9$, and where 
$t_n=n\Delta t$, $N = \left\lfloor \dfrac{T_{\text{learn}}}{\Delta t} \right\rfloor$;
  \item train a GP  using collected pairs of data,
    including optimization over hyperparameters such as lengthscale,
and set $\cgp$ to be the mean of this GP.
\end{enumerate}

Note that in step (c) we exploit the 
statistical invariance w.r.t. circular index 
shifting; this enables us to  
collect $K$ pairs in one time-step.
For numerical implementation of the GPR, we used the \texttt{scikit-learn}
package~\cite{scikit-learn}; the kernel of the GP was chosen as the sum 
of a standard radial basis function, and white noise kernels, setting 
the noise level in the latter to $0.5$; the length scale was optimized over.
Out of roughly $30000$ points obtained in step (c), we 
subsampled $500$ uniformly at random, without replacement, 
to train the GP.
The result of such a procedure is shown in Fig.~\ref{fig:gpr}.

\end{document}